\documentclass[12pt]{article}

\usepackage{graphicx}
\usepackage{subfig}
\usepackage{sidecap}
\usepackage{amsfonts}
\usepackage{amsmath}
\usepackage{amssymb}

\usepackage{psfrag}
\usepackage{float}

\usepackage{algorithm}
\usepackage{algpseudocode}
\usepackage{epstopdf}

\usepackage{soul}

\newtheorem{thm}{Theorem}[section] 
\newtheorem{lem}[thm]{Lemma} 
\newtheorem{cor}[thm]{Corollary}

\usepackage{bm}

\newcommand{\fraction}[2]{\textstyle\frac{#1}{#2}}
\newcommand{\Tr}{\textrm{Tr}}
\newcommand{\R}{\mathbb{R}}

\newcommand{\nh} {\widehat{\nabla}}

\newcommand{\bg}{\widehat{\nabla}F}
\newcommand{\Ss}{{\cal S}}
\newcommand{\Sh}{{\cal S}_H}

\newcommand{\bbh}{{b_H}}
\newcommand{\bb}{{b}}
\newcommand{\bM}{{M}}
\newcommand{\bL}{{L}}
\newcommand{\hti}{\mbox{$\tilde{m}$}}
\newcommand{\dtt}{\mbox{\rm det}}

\newcommand{\Bkb}{\mbox{$B_t^{(0)}$}}

\newcommand{\defeq}{\stackrel{\triangle}{=}}
\usepackage{soul,color}

\newcommand{\comment}[1]{}

\begin{document}
\title{ A Stochastic Quasi-Newton Method for Large-Scale Optimization}
\author{R. H. Byrd\thanks{Department of Computer Science, University of Colorado,
        Boulder, CO, USA.  This author was supported by National Science Foundation
        grant DMS-1216554 and Department of Energy grant DE-SC0001774.}
        \and 
        S.L. Hansen\thanks{Department of Engineering Sciences and Applied Mathematics, Northwestern University, 
       Evanston, IL, USA.  This author was supported by National Science Foundation grant DMS-0810213.} 
        \and   
       Jorge Nocedal \thanks{Department of Industrial Engineering and Management Sciences, Northwestern University, 
       Evanston, IL, USA.  This material is based upon work supported by the U.S. Department of Energy Office of Science, Office of Advanced
        Scientific Computing Research, Applied Mathematics program under Award Number FG02-87ER25047. This author was also supported by
        the Office of Naval Research award N000141410313.} 
       \and
        Y. Singer \thanks{Google Research}
      }

\maketitle
\begin{abstract}{ The question of how to incorporate curvature information in stochastic approximation methods is challenging. The direct application of classical quasi-Newton updating techniques for deterministic optimization leads to noisy curvature estimates that have harmful effects on the robustness of the iteration.  In this paper, we propose a stochastic quasi-Newton method that is efficient, robust and scalable. It employs the classical BFGS update formula in its limited memory form, and is based on the observation that it is beneficial to collect curvature information pointwise, and at regular intervals, through  (sub-sampled) Hessian-vector products. This technique differs from the classical approach that would compute differences of gradients at every iteration, and where controlling the quality of the curvature estimates can be difficult.  We present numerical results on problems arising in machine learning that suggest that the proposed method shows much promise. 
}
\end{abstract}
\newpage

\section{Introduction}
\label{formulation}
\setcounter{equation}{0}

 In many applications of machine learning, one constructs very large models from
massive amounts of training data. Learning such models 
imposes high computational and memory demands on the optimization algorithms
employed to learn the models. In some applications, a full-batch (sample
average approximation) approach is feasible and appropriate. However, in most large scale
learning problems, it is imperative to employ stochastic approximation
algorithms that update the prediction model based on a relatively small subset of
the training data.  These algorithms are particularly
suited for settings where data is perpetually streamed to the learning process;
examples include computer network traffic, web search, online advertisement,
and sensor networks.

The goal of this paper is to propose a quasi-Newton method that operates in
the stochastic approximation regime. We employ the well known
limited memory BFGS updating formula, and
show how to collect second-order information that is reliable enough to
produce stable and productive Hessian approximations. The key is to compute
average curvature estimates at regular intervals using (sub-sampled)
Hessian-vector products. This ensures sample uniformity and
avoids the potentially harmful effects of differencing noisy gradients.  

The problem under consideration is the minimization of a convex stochastic function, 
\begin{equation}  \label{ff}  
     \min_{w \in \R^n} F(w)   =  \mathbb{E} [f(w; \xi)]  ,
\end{equation}
where $\xi$ is a random variable. Although problem \eqref{ff} arises in other
settings, such as simulation optimization \cite{asmussen2007stochastic}, we assume for concreteness that $\xi$
is a random instance consisting of an input-output pair $(x,z)$. The vector
$x$ is typically referred to in machine learning as the input representation while $z$ as the
target output. In this setting, $f$ typically takes the  form

\begin{equation}   \label{pp}
  f(w; \xi)= f(w;x_i, z_i) = \ell (h(w; x_i); z_i),
\end{equation}
where $\ell$ is a loss function into $\R_+$, and $h$ is a
prediction model parametrized by $w$. The collection of input-output
pairs $\{(x_i, z_i)\}$, $i = 1,\cdots,N$ is   referred to as the
training set. The objective function \eqref{ff} is defined using the empirical
expectation 
\begin{equation} \label{empirical}
        F(w) = \frac{1}{N}\sum_{i=1}^N f(w;x_i, z_i).
\end{equation}

In learning applications with very large amounts of  training data, it is
common to use a mini-batch stochastic gradient based on $b \defeq |\Ss| \ll N$
input-output
instances, yielding the following estimate
\begin{equation}   \label{bat}
\bg(w) = \frac{1}{b} \sum_{i \in \Ss}  \nabla f(w;x_i, z_i) ~.
\end{equation}
The subset $\Ss$ $\subset \{ 1, 2, \cdots, N\}$ is randomly chosen, with 
$b $ sufficiently small so that the algorithm operates in the stochastic
approximation regime. Therefore, the stochastic estimates of the gradient are
substantially faster to compute than a gradient based on the entire  training set.

Our optimization 
method employs iterations of the form
\begin{equation}   \label{rm}
   w^{k+1} = w^k - \alpha^k B_k^{-1} \bg(w^k) ~,
\end{equation}
where $B_k$ is a symmetric positive definite approximation to the Hessian
matrix $\nabla^2 F(w)$, and $\alpha^k>0$. Since the stochastic gradient is not an accurate approximation to the gradient of \eqref{empirical} 
it is essential  (to guarantee convergence) that the steplength parameter $\alpha^k \rightarrow 0$. In our experiments and analysis,  $\alpha^k$ has the form $\alpha^k = \beta/k$, where $\beta >0$ is
given, but other choices can be employed.

A critical question is how to construct 
the Hessian approximation in a stable and efficient manner.
For the algorithm to be scalable, it must  update the inverse matrix $H_k= B_k^{-1}$ directly, so that \eqref{rm} can be implemented as
\begin{equation}   \label{direct}
        w^{k+1} = w^k - \alpha^k H_k \bg(w^k) .
 \end{equation}
 Furthermore, this step computation should require only $O(n)$ operations, as in limited memory quasi-Newton methods for deterministic optimization. 

If we set $H_k=I$ and $\alpha^k= \beta/k$ in \eqref{direct}, we recover the classical Robbins-Monro method \cite{RobMon51}, which is also called the \emph{stochastic gradient descent method}. Under standard convexity assumptions, the number of iterations needed by this method to compute an $\epsilon$-accurate solution is of order
$
      {n \nu\kappa^2}/{\epsilon} \, ,
 $ 
 where $\kappa$ is the condition number of the Hessian at the optimal solution, $\nabla^2 F(w^*)$, and $\nu$ is a parameter that depends on both the Hessian matrix and the gradient covariance matrix; see \cite{murata1998statistical,BottouBosq08}. Therefore, the stochastic gradient descent method is adversely affected by ill conditioning in the Hessian. In contrast, it is shown by Murata \cite{murata1998statistical} that setting $H_k=\nabla^2F(w^*)^{-1}$ in \eqref{direct}
completely removes the dependency on $\kappa$ from the complexity estimate. 
Although the choice $H_k=\nabla^2 F(w^*)^{-1}$ is not viable in practice, it suggests that an appropriate choice of $H_k$ may result in an algorithm that improves upon the stochastic gradient descent method.

In the next section, we present a stochastic quasi-Newton method of the form \eqref{direct} that is designed for large-scale applications. It employs the limited memory BFGS 
update, which is defined in terms of correction pairs $(s, y)$ that provide an estimate of the curvature of the objective function $F(w)$ along the most recently generated directions. We propose an efficient  way of defining these correction pairs that yields curvature estimates that are not corrupted by the effect 
of differencing the noise in the gradients. Our numerical experiments using problems arising in machine learning, suggest that the new method is robust and efficient.

The paper is organized into 6 sections. The new algorithm is presented in section~\ref{sec:qn}, and its convergence properties are discussed in section~\ref{analysis}. Numerical experiments that illustrate the practical performance of the algorithm are reported in section~\ref{numerical}. 
A literature survey on  related stochastic quasi-Newton methods is given in section~\ref{related}. The paper concludes in section~\ref{final} with some remarks about the contributions of the paper.
 
\bigskip\noindent\emph{Notation.}  The terms Robbins-Monro method, stochastic approximation (SA) method, and stochastic gradient descent (SGD) method are used in the literature to denote (essentially) the same algorithm. The first term is common in statistics, the second term is popular in the stochastic programming literature, and the acronym SGD is standard in machine learning. We will use the name stochastic gradient descent method (SGD) in the discussion that follows.

\section{A stochastic quasi-Newton method}  \label{sec:qn}
\setcounter{equation}{0}

The success of quasi-Newton methods for deterministic optimization lies in the fact that they construct curvature information during the course of the optimization process, and this information is good enough to endow the iteration with a superlinear rate of convergence. In the classical BFGS method \cite{Flet87} for minimizing a deterministic function $F(w)$, the new inverse approximation $H_{k+1}$ is uniquely determined by the previous approximation $H_k$ and the correction pairs
 \[
        y_k = \nabla F(w^{k+1}) - \nabla F(w^k), \quad s_k = w^{k+1} - w^k .
 \]
 Specifically,
\[
    H_{k+1} = (I- \rho_{k} s_k y_k^T) H_k (I- \rho_k  y_k s_k^T) + \rho_k s_k s_k^T \quad\mbox{with} \ \ 
    \rho_k = \frac{1}{y_k^T s_k}.
\]
This BFGS update is well defined as long as the curvature condition $y_k^T s_k >0$ is satisfied, which is always the case when  $F(w)$ is strictly convex.

For large scale applications, it is necessary to employ a limited memory variant that is scalable in the number of variables, but enjoys only a linear rate of convergence. This so-called L-BFGS method \cite{mybook} is considered generally superior to the steepest descent method for deterministic optimization: it produces well scaled and productive search directions that yield an approximate solution in fewer iterations and function evaluations. 


When extending the concept of limited memory quasi-Newton updating  to the stochastic approximation regime
 it is not advisable to mimic the classical approach for deterministic optimization and update the model based on information from only one iteration.
 This is because  quasi-Newton updating is inherently an overwriting process rather than an averaging process, and therefore the vector 
$y$ must reflect the action of the Hessian of the  entire objective $F$ given in \eqref{ff} --- something that is not achieved by differencing stochastic gradients \eqref{bat} based on small samples.  

 We propose that an effective approach to achieving stable Hessian approximation is to \emph{decouple} the stochastic gradient and curvature estimate calculations.
Doing so provides the opportunity to use a different sample 
subset for defining $y$ and the flexibility to add new  curvature estimates at regular intervals instead of at each iteration.
In order to emphasize that the curvature estimates are updated at a different schedule than the 
gradients, we use the subscript $t$ to denote the number of times a new $(s,y)$ pair has been calculated; this 
differs from the superscript $k$ which counts the number of gradient calculations and variables updates. 

 The displacement $s$ can be computed based on a collection of average iterates.
Assuming that new curvature estimates are calculated every $L$ iterations, we define $s_t$ as the difference of disjoint 
averages between the $2L$ most recent iterations: 
\begin{equation}
s_t = \bar{w}_t - \bar{w}_{t-1}, \quad\mbox{where} \quad \bar{w}_t = \sum_{i=k-L}^k w^i  ,\label{newcor}
\end{equation}
(and $ \bar{w}_{t-1}$ is defined similarly).
In order to avoid the potential harmful effects of gradient differencing when $\| s_t \|$ is small, we chose 
to compute $y_t$ via a Hessian vector product, 

\begin{equation}           
             y_t  =  \nh^2 F(\bar{w}_t)  s_t, \label{newcor2}
\end{equation}
i.e., by approximating differences in gradients via a first order Taylor expansion,
where $\nh^2 F$ is a sub-sampled Hessian defined as follows.
Let $\Sh \subset \{ 1, \cdots, N\}$ be a  randomly chosen subset of the training examples and let
\begin{equation}   \label{ssh}
  \nh^2F(w) \defeq
	\frac{1}{b_H} \sum_{i \in {\cal S}_{H}} \nabla^2 f(w;x_i,z_i) ,
\end{equation}
where $b_H$ is the cardinality of $\Sh$.

We  emphasize that the matrix $\nh^2 F(\bar w_t)$ is never
constructed explicitly when computing $ y_t$ in \eqref{newcor2}, rather, the
Hessian-vector product can be coded directly.  To provide useful curvature information, $\Sh$ should be relatively large (see section~\ref{numerical}),  regardless of the size of $b$.
%

The pseudocode of the complete method is given in
Algorithm~\ref{alg1}. 
\newpage

\bigskip

\begin{algorithm}
  \caption{Stochastic Quasi-Newton Method (SQN)}
  \label{alg1}
  {\bf Input:} initial parameters $w^1$, positive integers $M, L$, and
	step-length sequence $\alpha^k {>0}$   \\
  \begin{algorithmic}[1]
  \State Set $t = -1$ \Comment{Records number of correction pairs currently computed}
  \State $\bar{w}_t =  0$

  \For{ $k=1, \ldots, $}
     \State Choose a sample $\Ss$ $\subset \{ 1, 2, \cdots, N\}$
  	\State
	 Calculate stochastic gradient  $\bg(w^k)$ as defined in \eqref{bat}
	\State $\bar{w}_t =\bar{w}_t+ w^k$
	\If{ $k \leq 2L$ }
		\State $w^{k+1} = w^{k}-\alpha^k \bg(w^k)$ \Comment{Stochastic gradient iteration}
	\Else
		\State $w^{k+1} = w^{k}-\alpha^k H_{t} \bg(w^k)$, where $H_t$ is defined by Algorithm \ref{alg2}  \label{sqnDir} 
	\EndIf
	
	\If{mod$(k,L)=0$} \Comment{Compute correction pairs every $L$ iterations}
		\State $t=t+1$ 
		\State $\bar{w}_t = \bar{w}_t/L$
		\If{ $t > 0$} 
		\State Choose a sample $\Sh \subset \{ 1, \cdots, N\}$ to define $ \nh^2 F(\bar{w}_t) $ by  \eqref{ssh}
		\State Compute   \label{cpr}
		$      s_{t}=(\bar{w}_t-\bar{w}_{t-1}), \  y_{t}= \nh^2 F(\bar{w}_t) (\bar{w}_t-\bar{w}_{t-1}) $	\Comment{correction pairs}	
		\EndIf
		\State $\bar{w}_t=0$
	\EndIf
  \EndFor
  \end{algorithmic}
\end{algorithm}

\bigskip
 Using the averages \eqref{newcor} is not essential. One could also define $s$ to be the difference between two recent iterates.

The  L-BFGS step computation in  Step~\ref{sqnDir}  follows standard practice \cite{mybook}. Having chosen a memory parameter $M$, the matrix $H_t$ is defined as the result of applying $M$ BFGS updates to an initial matrix  using the $M$ most recent correction pairs 
$\{ s_j , y_j\}_{j=t-M+1}^t$ computed by the algorithm. This procedure is mathematically described as follows.
\newpage

\begin{algorithm}
  \caption{Hessian Updating}
  \label{alg2}
  {\bf Input:}  Updating counter $t$, memory parameter  $M$,  and  correction pairs $(s_j, y_j)$, 
  $j= t- \tilde{m}+1, \ldots t$, where $\hti = \min \{t,M\}$. \\
{\bf Output:} new matrix $H_t$
  \begin{algorithmic}[1]
   \State Set $
 H = \ ({s}_t^T {y}_t)/(y_t^T y_t) I, 
 $
where ${s}_t$ and  ${y}_t$ are computed in Step~\ref{cpr} of Algorithm~\ref{alg1}.

  \For{ $j=t- \tilde{m}+1,...,t$}
  	\State $ \rho_{j} = 1/y_j^Ts_j.$
	\State Apply BFGS formula:
	\begin{equation}  \label{lbfgs2}
    H \leftarrow (I- \rho_{j} s_j y_j^T) H (I- \rho_{j}  y_j  s_j^T) + \rho_j  s_j s_j^T \quad
\end{equation}

  \EndFor
    \State\Return $H_t \leftarrow H$
  \end{algorithmic}
\end{algorithm}

\bigskip

In practice, the quasi-Newton matrix $H_t$ is not formed explicitly; to compute the product 
$H_{t} \bg(w^k)$  in Step 10 of Algorithm~\ref{alg1} one employs a formula based on the structure of the 2-rank BFGS update. This formula, commonly called the two-loop recursion,  computes the step directly from the correction pairs and stochastic gradient as described in \cite[\S 7.2]{mybook}.  


In summary, the algorithm builds upon the strengths of BFGS updating, but deviates from the classical method in that the correction pairs 
$(s, y)$ are based on  \emph{sub-sampled} Hessian-vector products computed at regularly spaced intervals, which amortize their cost. Our task in the remainder of the paper is  to argue that even with the extra computational cost of Hessian-vector products \eqref{newcor2} and the extra cost of computing the iteration \eqref{direct}, the stochastic quasi-Newton method is competitive with the SGD method in terms of computing time (even in the early stages of the optimization), and is able to find a lower objective value.

\subsection{Computational Cost}\label{stochb1}

Let us compare the cost of the stochastic gradient descent method
\begin{equation}   \label{sgd}
     w^{k+1}= w^k - \frac{\beta}{k} \bg(w^k) \ \ \ \ \ \textrm{(SGD)}
\end{equation}
and the stochastic quasi-Newton method
\begin{equation}   \label{stqn}
     w^{k+1}= w^k - \frac{\beta}{k} H_t \bg(w^k) \ \ \ \textrm{(SQN)}
\end{equation}
given by Algorithm~\ref{alg1}.

The quasi-Newton matrix-vector product in \eqref{stqn} requires approximately $4Mn$ operations \cite{mybook}. To measure the cost of the gradient and Hessian-vector computations, let us consider one particular but representative example, namely the binary classification test problem tested in section~\ref{numerical}; see \eqref{ysi}. In this case, the component function $f$ in \eqref{pp} is given by
\[ 
    f(w; x_i,z_i) =  z_i\log(c(w;x_i)) + (1-z_i)\log\left(1-c(w; x_i)\right)
\]
where
\begin{equation} \label{blc}
    c(w;x_i)=\frac{1}{1+\exp(-x_i^Tw)}, \quad x_i\in\mathbb{R}^{n}, \ w\in\mathbb{R}^{n}, \ z_i\in\left\{0,1\right\}.
\end{equation}
The gradient and Hessian-vector product of $f$ are given by,
\begin{align}
\nabla f(w; x_i,z_i) &= (c(w;x_i) - z_i) x_i    \label{five1} \\
\nabla^2 f(w; x_i,z_i) s &= c(w; x_i) (1-(c(w; x_i))( x_i^T  s) x_i . \label{five2}
\end{align}
The evaluation of the function $c(w;x_i)$ requires approximately $n$ operations (where we follow the convention of counting a multiplication and an addition as an operation). 
Therefore, by \eqref{five1} the cost of evaluating one batch gradient is approximately $2 \bb n$, and the cost of computing the Hessian-vector product $\widehat \nabla^2 F(\bar{w}_t) s_t$ is  about $3\bbh n$. 
This assumes these two vectors are computed independently. If the Hessian is computed at the same point where we compute a gradient and $\bb \geq \bbh$ then $c(w;x_i)$ can be reused for a savings of $\bbh n$.

Therefore, for binary logistic  problems the total number of floating point operations of the stochastic quasi-Newton iteration \eqref{stqn} is approximately 
 \begin{equation}   \label{counts}
          2\bb n  + 4\bM n +  3 \bbh n/\bL .
 \end{equation}
On the other hand, the cost associated with the computation of the SGD step is only $b n$. At first glance it may appear that the SQN method is prohibitively expensive, but this is not the case 
when using the values for  $\bb$, $\bbh$, $\bL$ and $\bM$ suggested in this paper.
To see this, note that
 \begin{equation}    \label{cratios}
     \frac{\mbox{cost of SQN iteration}}{\mbox{cost of SGD iteration}} =  1 + \frac{2\bM}{\bb}+ \frac{ 2 \bbh}{3 \bb  \bL}.~
\end{equation}
In the experiments reported below, we use $\bM=5$, $\bb=50, 100, \ldots,$ $\bL=10$ or 20, and choose $\bbh \geq 300$. For such parameter settings, the additional cost of the SQN iteration is small relative to the cost of the SGD method.

For the multi-class logistic regression problem described in section~\ref{speechs}, the costs of gradient and Hessian-vector products are slightly different. Nevertheless, the relative cost of the SQN and SGD iterations is similar to that given in \eqref{cratios}.

The quasi-Newton method can take advantage of parallelism. Instead of employing the two-loop recursion mentioned above 
to implement the limited memory BFGS step computation in step \ref{sqnDir} of Algorithm~\ref{alg1}, we can employ the compact form of limited memory BFGS updating \cite[\S 7.2]{mybook} in which $H_t$ is represented as the outer product of two matrices. This computation can be parallelized and its effective cost is around $3n$ operations, which is smaller than the $4 \bM n$ operations assumed above. The precise cost of parallelizing the compact form computation depends on the computer architecture, but is in any case independent of $\bM$.
%

Additionally, the Hessian-vector products can be computed in parallel with the main iteration \eqref{stqn}  if we allow freedom in the choice of the point $\bar w_t$ where \eqref{newcor2} is computed.
The choice of this point is not delicate since it suffices to estimate the average curvature of the problem around the current iterate, and hence the computation of \eqref{newcor2} can lag behind the main iteration.
In such a parallel setting, the computational overhead of Hessian-vector products may be negligible. 

The SQN method contains several parameters, and we provide the following guidelines on how to select them. First, the minibatch size $b$ is often dictated by the experimental set-up or the computing environment, and we view it as an exogenous parameter. A key requirement in the design of our algorithm is that it should  work well with any value of $b$.  Given $b$, our experiments show that it is most efficient if the per-iteration cost
of updating, namely $b_H/L$,  is less than the cost of the stochastic gradient $b$, with the ratio $Lb/b_H$ in the range [2,20].   
The choice of the  parameter  $M$ in L-BFGS updating is similar as in deterministic optimization; the best value is problem dependent but values in the  range [4,20] are commonly used. 

\section{Convergence Analysis}  \label{analysis}
\setcounter{equation}{0}

In this section,  we analyze the convergence properties of the stochastic quasi-Newton method. We assume that the objective function $F$  is strongly convex and twice continuously differentiable. The first assumption may appear to be unduly strong because in certain settings (such as logistic regression) the component functions $f(w;x_i,z_i)$ in \eqref{empirical} are  convex, but not strongly convex. However, since the lack of strong convexity can lead to  very slow convergence, it is common in practice to either add an $\ell_2$ regularization term, or choose the initial point (or employ some other mechanism) to ensure that the iterates remain in a region where the $F$ is strongly convex.   If regularization is used, the objective function  takes the form
\begin{equation} \label{remp}
         \fraction{1}{2} \sigma \|w\|^2 + \frac{1}{N}\sum \limits_{i=1}^N f(w;x_i, z_i),   \quad\quad\mbox{with} \ \  \sigma>0,
\end{equation} 
and the sampled Hessian \eqref{ssh} is  
\begin{equation}   \label{ssh2}
    \sigma I +
	\frac{1}{b_H} \sum_{i \in {\cal S}_{H}} \nabla^2 f(w;x_i,z_i) .
\end{equation}
In this paper, we do not specify the precise mechanism by which the strong convexity is ensured.
The  assumptions made in our analysis  are  summarized as follows. 

\bigskip\noindent
\textbf{Assumptions 1}

\noindent\emph{
  (1) The objective function $F$ is twice continuously
           differentiable.}
           
            \smallskip     
 \noindent \emph{ (2) There exist positive constants $\lambda $ and $\Lambda$ such that
 \begin{equation}   
          \lambda  I \prec \nh^2F(w) \prec \Lambda I , \label{unif}
\end{equation}
            for all $w \in \mathbb{R}^n$, and all ${\cal S}_{H}  \subseteq \{1,\cdots, N\}$. (Recall that ${\cal S}_{H}$ appears in the
  definition \eqref{ssh} of \ $ \nh^2F(w)$.)
            }

\bigskip\noindent
These assumptions imply that the entire Hessian $\nabla^2 F(w)$ satisfies \eqref{unif} for all $w \in \mathbb{R}^n$, and that $F$ has a unique   minimizer $w^*$. 
If the matrices  $\nabla^2 f(w;,x_i, z_i)$ are nonnegative definite and uniformly bounded and 
we implement $\ell_2$ regularization as in (\ref{remp})--(\ref{ssh2}) then part (2) of Assumptions~1 is satisfied.    

We first show that   the Hessian approximations  generated by the SQN method have eigenvalues that are uniformly bounded above and away from zero. 
%

\begin{lem}  \label{lem1}
If Assumptions~1 hold, 
 there exist constants $0 < \mu_1\leq  \mu_2$ such that the Hessian approximations $\{H_t\}$ generated by  Algorithm~\ref{alg1} satisfy 
\begin{equation}
   \mu_1 I \prec H_t \prec \mu_2 I  , \quad\quad\mbox{ for} \ \ t=1,2, \ldots 
   \label{linc}
\end{equation}
\end{lem}

\medskip\noindent{\bf Proof:} 
Instead of analyzing the inverse Hessian approximation $H_k$, we will study the direct Hessian approximation $B_k$ (see \eqref{rm}) because this allows us to easily quote a result from the literature. In this case, the limited memory quasi-Newton updating formula is given as follows: 

\medskip i)  
Set $B_t^{(0)}= \frac{ y_t^T  y_t}{ {s}_t^T {y}_t} I $, and  $\tilde m = \min \{t, M\}$;

ii) for  $i=0,..,\tilde m-1 $ set $j=t-\tilde{m} +1+i$ and compute
\begin{equation}
    B_{t}^{(i+1)} = B_t^{(i)} - \frac{B_t^{(i)} s_{j} s_{j} ^T B_t^{(i)} }{ s_{j} ^T B_t^{(i)} s_{j} }
                     + \frac{y_{j}  y_{j} ^T}{y_{j} ^Ts_{j} }.
                                                          \label{lbfgs}
\end{equation}

(iii)
Set $B_{t+1}=B_t^{(\tilde m)}$.
\medskip

\noindent By \eqref{newcor2}
\begin{equation}    \label{lime}
        s_j = \bar w_j  - \bar w_{j-1} ,  \qquad  y_j = \nh^2 F(\bar w_j )  s_j  .  
\end{equation}
and thus by \eqref{unif} 
\begin{equation}
    \lambda \|s_j\|^2 \leq y_j^Ts_j \leq  \Lambda \|s_j\|^2.   \label{r1}
\end{equation}
Now, 
\[
   \frac{\|y_j\|^2}{y_j^Ts_j} = \frac{s_j^T  {\nh^2 F(\bar w_t )^2} s_j} 
   {s_j^T \nh^2 F(\bar w_t ) s_j} ,        
\]
and since $ \nh^2 F(\bar w_t )$ is symmetric and positive definite, it has a square root so that
\begin{equation}
   \lambda \leq \frac{\|y_j\|^2}{y_j^Ts_j}  \leq \Lambda.           \label{r2}
\end{equation}
This proves that the eigenvalues of the matrices $B_t^{(0)}= \frac{ y_t^T  y_t}{ {s}_t^T {y}_t} I$ at the start of the L-BFGS update cycles are bounded above and away from zero, for all $t$.

Let $\Tr(\cdot)$ denote the trace of a matrix. Then from (\ref{lbfgs}), (\ref{r2}) and the boundedness
of $\{\|\Bkb\|\}$, and setting  $j_i=t- \tilde m +i$,
\begin{eqnarray}
    \Tr(B_{t+1}) & \leq & \Tr(\Bkb) + \sum_{i=1}^{\tilde m}\frac{\|y_{j_i}\|^2}{y_{j_i}^T s_{j_i}}  \nonumber   \\
                 & \leq & \Tr(\Bkb) + \tilde{m} \Lambda                                        \nonumber \\
                 & \leq & M_3,                             \label{ltrace}
\end{eqnarray}
for some positive constant $M_3$. This implies that the largest eigenvalue of all matrices $B_t$ is bounded uniformly.

We now derive an expression for the determinant of $B_t$. It is shown by  Powell  \cite{Powe766} that
\begin{eqnarray}
   \dtt(B_{t+1}) & = & \dtt(\Bkb) \prod_{i=1}^{\tilde m} \frac{y_{j_i}^Ts_{j_i}}{s_{j_i}^T B_t^{(i-1)}s_{j_i}}  \nonumber \\
                 & = & \dtt(\Bkb) \prod_{i=1}^{\tilde{m}} \frac{y_{j_i}^Ts_{j_i}}{s_{j_i}^Ts_{j_i}} 
                        \frac{s_{j_i}^Ts_{j_i}}{s_{j_i}^TB_t^{(i-1)}s_{j_i}}.  \label{det1}
\end{eqnarray}
Since by (\ref{ltrace}) the largest eigenvalue of $B_t^{(i)}$ is  less than $M_3$, we have,
using  (\ref{r1}) and the fact that the smallest eigenvalue of $\Bkb$ is bounded away from zero,
\begin{eqnarray}
   \dtt(B_{t+1}) & \geq & \dtt(\Bkb)\left(\frac{\lambda}{M_3}\right)^{\tilde{m}}                 \nonumber  \\
                 & \geq & M_4,                               \label{det2} 
\end{eqnarray}
for some positive constant $M_4$. This shows that the smallest eigenvalue of the matrices $B_t$ is bounded away from zero, uniformly. Therefore, condition \eqref{linc} is satisfied. \hspace*{\fill}$\Box$\medskip

 
 Our next task is to establish global convergence. Rather than proving this result just for our SQN method, we analyze a more general iteration that covers it as a special case. We do so because the more general result is of interest beyond this paper and we  are unaware of a self-contained proof of this result in the literature (c.f.  \cite{sunehag2009variable}). 

We consider the  Newton-like iteration
\begin{equation}  \label{iteration}
w^{k+1}=w^k-\alpha^k H_k \nabla f(w^k, \xi^k),
\end{equation}
when applied to a strongly convex objective function $F(w)$. (As in \eqref{pp},  we used the notation $\xi = (x,z)$.) We assume that  the eigenvalues of $\{H_k\}$ are uniformly bounded above and away from zero, and that $E_{\xi^k}[\nabla f(w^k, \xi^k)]= F(w^k)$.  Clearly Algorithm~\ref{alg1} is a special case of \eqref{iteration} in which $H_k$ is constant for $L$ iterations. 

We make the following assumptions about the  functions $f$ and $F$.


\bigskip\noindent
\textbf{Assumptions 2}

\noindent\emph{
(1) $F(w)$ is twice  continuously differentiable. }

\smallskip \noindent \emph{(2) There  exist positive constants $\lambda $ and $\Lambda$ such that, for all $w \in \mathbb{R}^n,$
\begin{equation}
          \lambda I \prec \nabla^2F(w) \prec \Lambda I . \label{uniff}
\end{equation} }
\noindent \emph{(3) There is a constant $\gamma$ such that, for all $w \in \mathbb{R}^n$,
           \begin{equation}  \label{varo}
E_{\xi}[ \| \nabla f(w^k, \xi)\|)]^2  \leq \gamma^2 .  
\end{equation}
}

\medskip\noindent
For convenience, we define $\alpha^k = \beta /k$, for an appropriate choice of $\beta$, rather than assuming well known and more general conditions $\sum \alpha^k=\infty$, $\sum (\alpha^k)^2 <\infty$. This allows us to provide a short proof  similar to the analysis of Nemirovski et al. \cite{nemirovski2009robust}. 
\begin{thm}  \label{thm1}
Suppose that Assumptions~2 hold. Let $w^k$ be the iterates generated by the Newton-like method \eqref{iteration}, 
 where for  $k=1,2, \ldots$
\begin{equation}    \label{mbounds}
      \mu_1 I \prec H_k \prec \mu_2 I ,  \qquad 0< \mu_1 \leq \mu_2 ,
\end{equation}
and
\[
    \alpha^k = \beta/{k} \quad \mbox{ with}\quad  \beta > 1/({2 \mu_1 \lambda}).  
\] 
Then for all $k \geq 1$,
\begin{equation}   \label{mainr}
    E[F(w^k)-F(w^\ast)] \leq Q(\beta)/k ,
\end{equation}
 where
\begin{equation}   \label{qbound}
     Q(\beta) =\max \left\{ {\Lambda \mu_2^2\beta^2 \gamma^2 \over 2(2 \mu_1 \lambda \beta -1)}, F(w^1) -F(w^\ast) \right\} .    
 \end{equation}
\end{thm}
\textbf{Proof.} We have that
\begin{eqnarray*}
F(w^{k+1})&= &F(w^k -\alpha^k H_k \nabla f(w^k, \xi^k))   \\
          &\leq &F(w^k) +\nabla F(w^k)^T(-\alpha^k H_k \nabla f(w^k, \xi^k)) + \fraction{\Lambda}{2}  \|\alpha^k H_k \nabla f(w^k, \xi^k)\|^2 \\
        &\leq & F(w^k) - \alpha^k \nabla F(w^k)^T H_k \nabla f(w^k, \xi^k)) +\fraction{\Lambda}{2}  (\alpha^k \mu_2 \| \nabla f(w^k, \xi^k)\|)^2  .
\end{eqnarray*}
Taking the expectation over all possible values of $\xi^k$ and recalling \eqref{varo} gives
\begin{eqnarray}
E_{\xi^k}[F(w^{k+1})] & \leq & F(w^k)  - \alpha^k \nabla F(w^k)^T H_k \nabla F(w^k) + \fraction{\Lambda}{2}  (\alpha^k \mu_2)^2 E_{\xi^k}[ \| \nabla f(w^k, \xi^k)\|]^2    \nonumber \\
& \leq&  F(w^k)  - \alpha^k \mu_1 \| \nabla F(w^k)\|^2 + \fraction{\Lambda}{2}  (\alpha^k \mu_2)^2  \gamma^2.   \label{firstexpect}
\end{eqnarray}
Now, to relate $F(w^k) -F(w^\ast)$ and $\| \nabla F(w^k)\|^2$, we use the lower bound in \eqref{uniff} to construct a minorizing quadratic for $F$ at $w^k$. For any vector $v \in \mathbb{R}^n$, we have
\begin{eqnarray}
F(v) &\geq &F(w^k) +\nabla F(w^k)^T (v-w^k) + \fraction{\lambda}{ 2} ||v-w^k||^2  \nonumber \\
     &\geq &F(w^k) +\nabla F(w^k)^T (-\fraction{1}{ \lambda} \nabla F(w^k) ) + {\lambda \over 2} ||{1\over \lambda} \nabla F(w^k)||^2  \nonumber \\
     &\geq &F(w^k)  - \fraction{1}{2 \lambda} || \nabla F(w^k)||^2  , \label{uber}
\end{eqnarray}
where the second inequality follows from the fact that $\hat v = w^k - \frac{1}{\lambda} \nabla F(w^k)$ minimizes the quadratic
$     q_k(v)= F(w^k) +\nabla F(w^k)^T (v-w^k) + \fraction{\lambda}{ 2} ||v-w^k||^2 .$
Setting $v=w^\ast$ in \eqref{uber} yields
\[
 2 \lambda [F(w^k) -F(w^\ast)] \leq \|\nabla F(w^k) \|^2 ,
\]
which together with  (\ref{firstexpect}) yields 
\begin{equation}   \label{wifi}
E_{\xi^k}[F(w^{k+1}) -F(w^\ast)] \leq   F(w^k)-F(w^\ast)  - 2 \alpha^k  \mu_1 \lambda [F(w^k -F(w^\ast)] + 
    \fraction{\Lambda}{2}  (\alpha^k \mu_2)^2 \gamma^2  .
\end{equation}
Let us define $\phi_k$ to be the expectation of $F(w^k) - F(w^\ast)$ over all choices 
 $\{ \xi^1,\xi^2, \ldots, \xi^{k-1}\}$ starting at $w^1$, which we write as 
\begin{equation}  \label{fullex}
       \phi_k = E[F(w^k) - F(w^\ast)].
\end{equation}
 Then equation \eqref{wifi} yields 
 \begin{equation} \label{fi3}
 \phi_{k+1} \leq  (1 - 2 \alpha^k  \mu_1 \lambda ) \phi_k + \fraction{\Lambda}{2}  (\alpha^k \mu_2)^2  \gamma^2 .
\end{equation}

We prove  the desired result \eqref{mainr} by induction. 
The result clearly holds for $k=1$. Assuming it holds for some value of $k$,
inequality (\ref{fi3}), definition \eqref{qbound}, and the choice of $\alpha^k$ imply
\begin{eqnarray*}
 \phi_{k+1} &\leq&  \left(1 -  \frac{2 \beta \mu_1 \lambda }{ k} \right) {Q(\beta) \over k} + {\Lambda\mu_2^2 \beta^2 \gamma^2 \over 2k^2 }   \\
  &=& {(k - 2 \beta  \mu_1 \lambda) Q(\beta) \over k^2} + {\Lambda\mu_2^2 \beta^2 \gamma^2 \over 2k^2 } \\
  &=& {(k - 1) Q(\beta) \over k^2} - {2 \beta  \mu_1 \lambda -1 \over k^2} Q(\beta) + {\Lambda\mu_2^2 \beta^2 \gamma^2 \over 2k^2 }\\
  &\leq& { Q(\beta) \over k+1} .
\end{eqnarray*} 
 \hspace*{\fill}$\Box$\medskip


%
\noindent We can now establish the convergence result.
\begin{cor} 
Suppose that Assumptions 1 and the bound \eqref{varo} hold. Let $\{w^k\}$ be the iterates generated by Algorithm \ref{alg1}.
Then there is a constant $\mu_1$ such that 
     $ \|H_k^{-1} \| \leq {1 / \mu_1 }$ ,
and if the steplength is chosen by
\[
    \alpha^k = \beta/{k} \quad \mbox{ where}\quad  \beta > {1}/{(2 \mu_1 \lambda)}, \quad \forall k ,
\] 
it follows that
\begin{equation}   
    E[F(w^k)-F(w^\ast)] \leq Q(\beta)/k ,
\end{equation}
for all $k$, where
\begin{equation}   
     Q(\beta) =\max \left\{ {\Lambda \mu_2^2\beta^2 \gamma^2 \over 2(2\mu_1 \lambda \beta -1)}, F(w^1) -F(w^\ast) \right\} .
 \end{equation}
\end{cor}
{\bf Proof:} Lemma~\ref{lem1} ensures that the Hessian approximation satisfies \eqref{linc}. Now, the iteration in Step~10 of Algorithm~\ref{alg1} is a special case of iteration \eqref{iteration}. Therefore, the result follows from Theorem~\ref{thm1}.
 \hspace*{\fill}$\Box$\medskip

\section{Numerical Experiments}  \label{numerical}
\setcounter{equation}{0}

In this section, we compare the performance of the stochastic gradient descent (SGD) method \eqref{sgd} and the stochastic quasi-Newton (SQN) method \eqref{stqn}
on three test problems of the form \eqref{pp}-\eqref{empirical} arising in supervised machine learning. The parameter $\beta>0$ is fixed at the beginning of each run, as discussed below, and the SQN method is implemented as described in Algorithm~\ref{alg1}.  

It is well known amongst the optimization and machine learning communities that the SGD method can be improved by choosing the parameter $\beta$ via a set of problem dependent heuristics \cite{plakhov2004stochastic, yousefian2012stochastic}. In some cases, $\beta_k$ (rather than $\beta$) is made to vary during the course of the iteration, and could even be chosen so that $\beta_k/k$ is constant, in which case only convergence to a neighborhood of the solution is guaranteed \cite{nedic2001convergence}. There is, however, no generally accepted rule for choosing $\beta_k$, so our testing approach is to consider the simple strategy of selecting the (constant) $\beta$ so as to give good performance for each problem. 

Specifically, in the experiments reported below, we tried several values for $\beta$ in \eqref{sgd} and \eqref{stqn} and chose a value for which increasing or decreasing it by a fixed increment results in inferior performance. This allows us to observe the effect of the quasi-Newton Hessian approximation $H_k$ in a controlled setting,  without the clutter introduced by elaborate step length strategies for $\beta_k$. 

\bigskip\noindent In the figures provided below, we use the following notation.

\begin{enumerate}
\item $n$: the number of variables in the optimization problem; i.e., $w \in \R^n$.
\item $N$: the number of training points in the dataset.
\item $\bb$: size of the batch used in the computation of the stochastic gradient $\widehat \nabla F(w)$ defined in \eqref{bat}; i.e., $\bb = | {\cal S}|$.
\item $\bbh$:  size of the batch used in the computation of Hessian-vector products \eqref{newcor2} and \eqref{ssh}; i.e., $\bbh = | {\cal S}_H|$.

\item $\bL$:   controls the frequency of limited memory BFGS updating. Every $\bL$ iterations a new curvature pair (\ref{newcor2}) is formed and the oldest pair is removed.
                                             
\item $\bM$: memory used in limited memory BFGS updating.
 
\item adp: Accessed data points. At each iteration the
SGD method evaluates the stochastic gradient $\nh F(w^k)$ using $\bb$ randomly chosen training points $(x_i, z_i)$,
so we say that the iteration accessed $\bb$ data points. On the other hand, an iteration of the stochastic BFGS method accesses $\bb + \bbh/\bL$ points.
 \item iteration: In some graphs we compare SGD and SQN iteration by iteration  (in addition to comparing them in terms of accessed data points).
 \item epoch:  One complete pass through the dataset.                                         
  \end{enumerate}

In our experiments, the stochastic gradient \eqref{bat} is formed by randomly choosing $\bb$ training points from the dataset without replacement. This process is repeated every epoch, which guarantees that all training points are equally used when forming the stochastic gradient. 
Independent of the stochastic gradients, the Hessian-vector products are formed by randomly choosing $\bbh$ training points from the dataset without replacement.


 \subsection{Experiments with Synthetic Datasets} \label{ysi}
 We first test our algorithm on a binary classification problem.
The objective function is given by
\begin{equation}   \label{rcvf}
    F(w) = -\frac{1}{N} \sum_{i=1}^N z_i\log(c(w;x_i)) + (1-z_i)\log\left(1-c(w; x_i)\right)  ,
\end{equation} 
where $c(w;x_i)$ is defined in (\ref{blc}).

The training points were generated randomly as described in \cite{41341}, with $N=7000$ and $n=50$. To establish a reference benchmark with a well known algorithm, we used the particular implementation \cite{41341} of one of the coordinate descent (CD) methods of Tseng and Yun~\cite{tseng2009coordinate}. 

\begin{figure}[H] 
 \centering
\includegraphics[scale=.65]{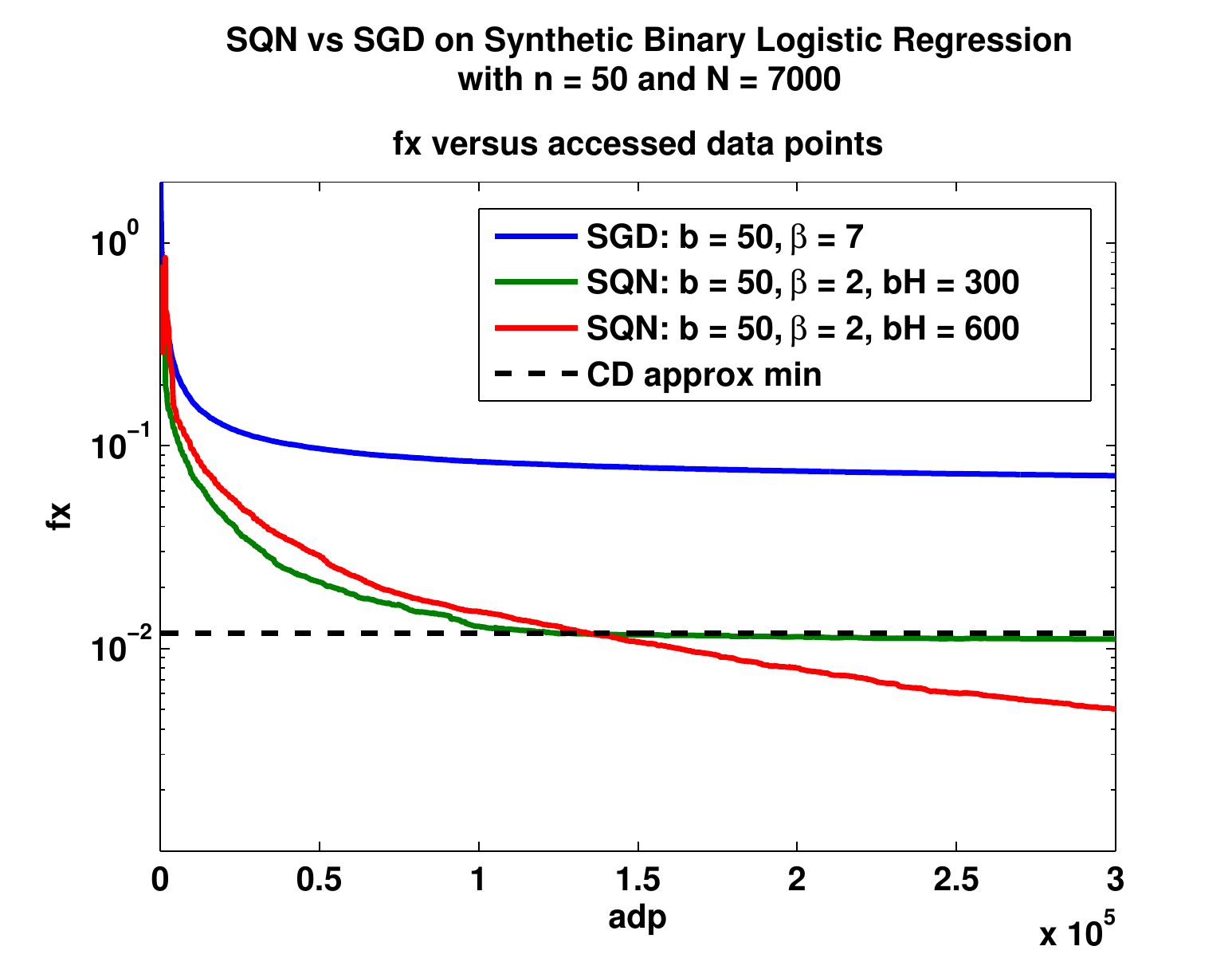}
\caption{Illustration of SQN and SGD on the synthetic dataset. The dotted black line  marks the best function value obtained by the coordinate descent (CD) method.
For SQN we set $\bM=10$, $\bL=10$ and $\bbh=300$ or 600.}
\label{random1}
\end{figure}

Figure~\ref{random1} reports the performance of SGD (with $\beta=7$) and SQN (with $\beta=2$),  as measured by accessed data points. Both methods use a gradient  batch size of $ \bb=50$; for SQN we display results for two values of the Hessian batch size $\bbh$, and set  $\bM=10$ and $\bL=10$. The vertical axis, labeled {\tt fx}, measures the value of the objective \eqref{rcvf}; the dotted black line marks the best function value obtained by the coordinate descent (CD) method mentioned above.
We observe that the SQN method with $\bbh = 300$ and 600 outperforms SGD,
and obtains the same or better objective value than the coordinate descent method. 

\begin{figure}[H] 
 \centering
\includegraphics[scale=.60]{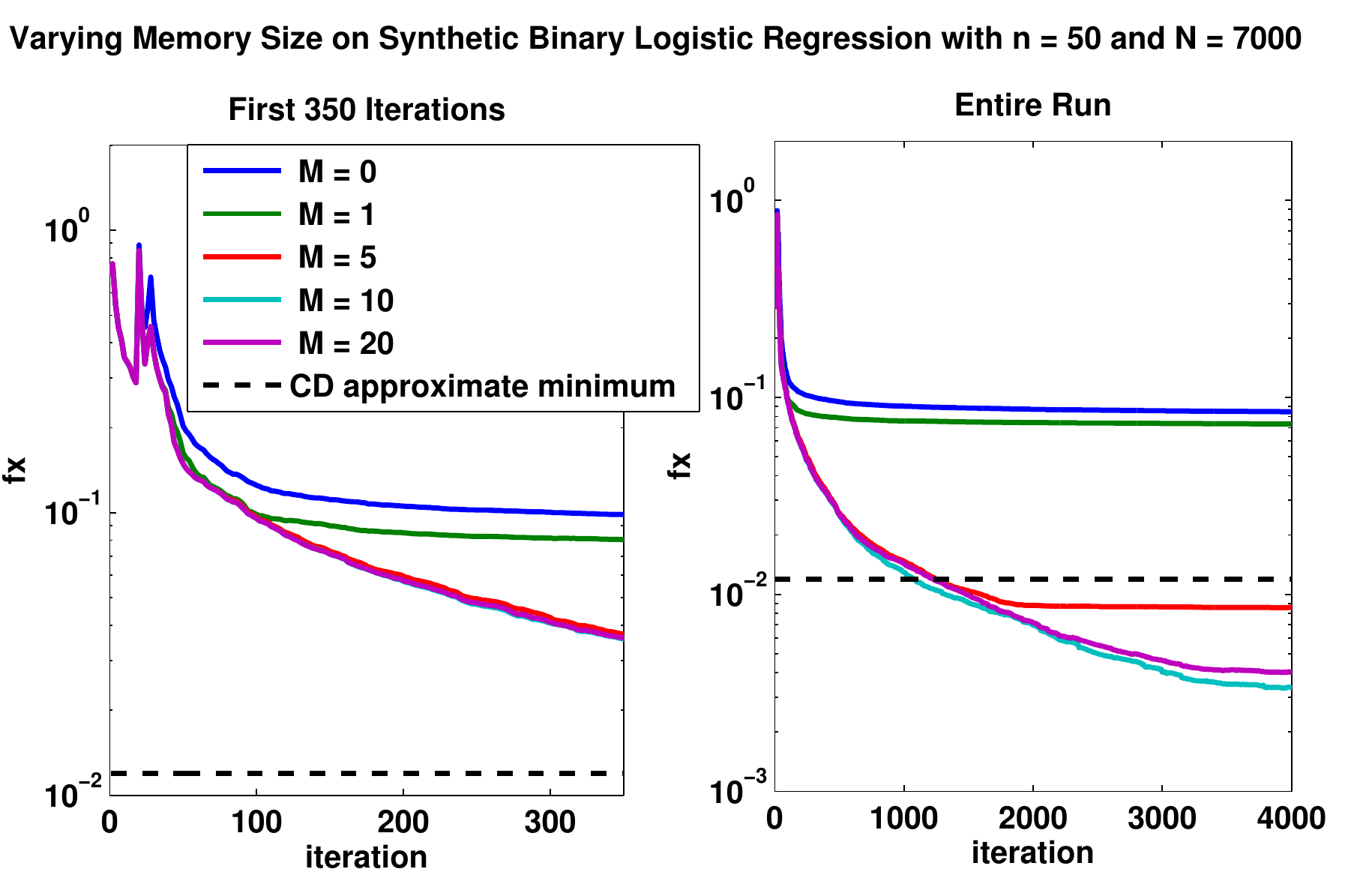}
\caption{Effect of the memory size $\bM$ in the SQN method. The figure on the left reports the first 4 epochs, while the figure on the right lets the algorithm run for more than 70 epochs to observe if the beneficial effect of increasing $\bM$ is sustained. Parameters settings are $\bb=50$, $\bbh=600$, and $\bL=10$.}
\label{random2}
\end{figure}

In Figure \ref{random2} we explore the effect of the  memory size $\bM$. 
Increasing $\bM$ beyond 1 and 2 steadily improves the performance of the SQN algorithm, both during the first few epochs (left figure), and after letting the algorithm run for many epochs (right figure). For this problem, a large memory size is
helpful in the later stages of the run.

\subsection{RCV1 Data Set}\label{rcv1s}

The RCV1 dataset \cite{lewis2004rcv1} is a composition of newswire articles produced by Reuters from 1996-1997. Each article was manually labeled into 4 different classes: Corporate/Industrial, Economics, Government/Social, and Markets. For the purpose of classification, each article was then converted into a boolean feature vector with a 1 representing the appearance of a given word. Post word stemming, this gives each feature vector a dimension of $n=112919$.

Each data point $x_i\in[0,1]^n$ is extremely sparse, with an average of 91 ($.013\%)$ nonzero elements. There are $N=688329$ training points. 
We consider the binary classification problem of predicting
whether or not an article is in the fourth class, Markets, and accordingly we have labels $z_i\in\left\{ 0, 1\right\}$.
We use logistic regression to model the problem, and define the objective function by equation (\ref{rcvf}).

In our numerical experiments with this problem, we used gradient batch sizes of $\bb = 50$, 300 or 1000, which respectively comprise $.0073\%$, $.044\%$ and $.145\%$ of the dataset. The frequency of quasi-Newton updates was set to $\bL=20$, a value that balances the aims of quickly retrieving curvature information and minimizing computational costs. For the SGD method we chose $\beta=5$ when $\bb=50$, and  $\beta=10$ when $\bb = 300$ or 1000; for the SQN method \eqref{stqn} we chose  $\beta=1$ when $\bb=50$, and $\beta=5$ when $\bb = 300$ or 1000, and we set $\bbh=1000$.

\begin{figure}[H]
 \centering
\includegraphics[scale=.55]{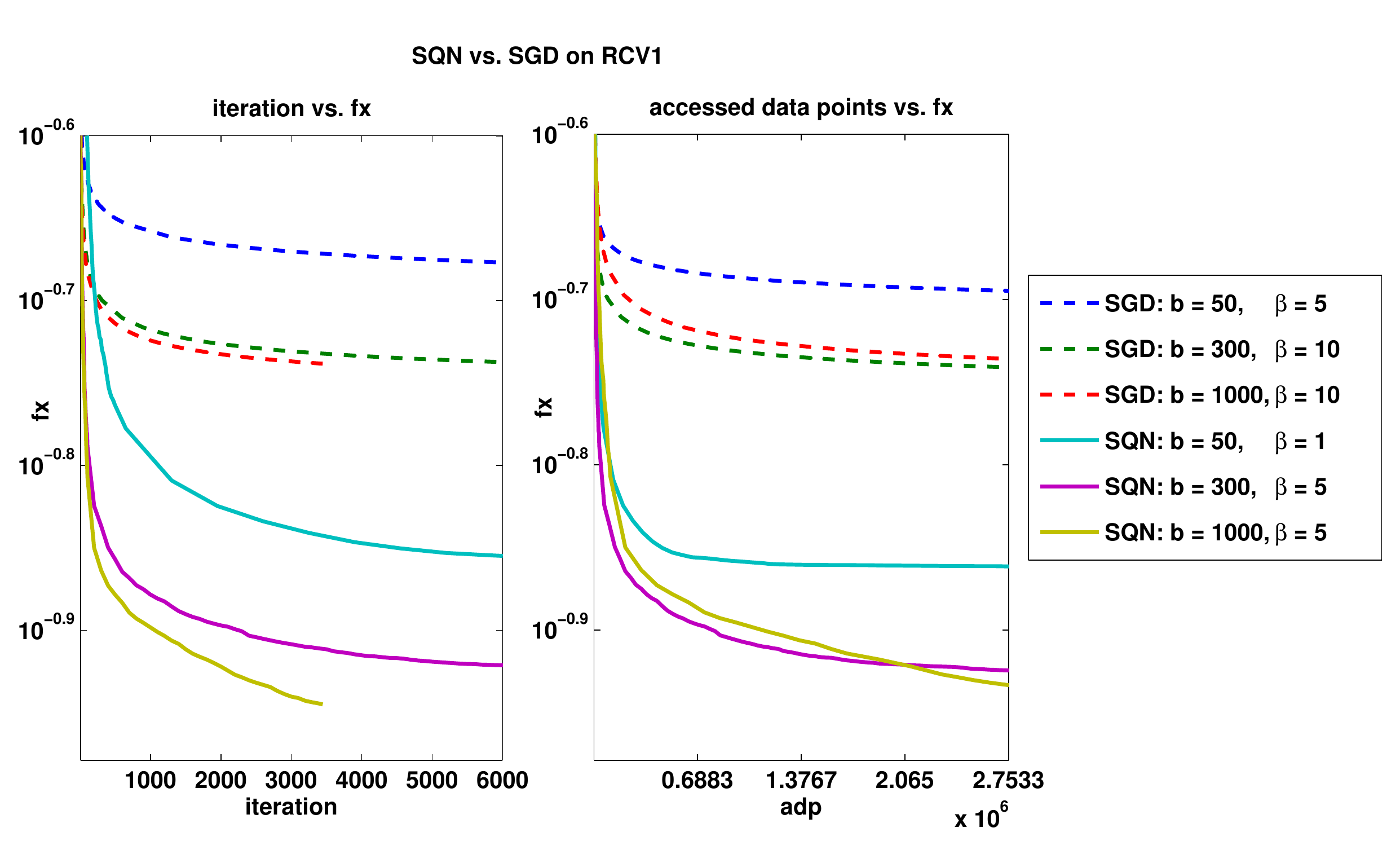}
\caption{Illustration on RCV1 problem. For SGD and SQN, $\bb$ is set to either 50, 300 or 1000, and for SQN we use $\bbh=1000$,  $\bM=5$, and $\bL=20$. The figures report training error as a function of iteration count or accessed data points. In the rightmost graph the tick marks on the x-axis (at 0.6882, 1.3767, $\ldots$) denote the epochs of SGD. }
\label{fig:rcv1sgdVsqn} 
\end{figure}

Figure~\ref{fig:rcv1sgdVsqn} reports the performance of  the two methods as measured by either iteration count or accessed data points. As before, the vertical axis, labeled {\tt fx}, measures the value of the objective \eqref{rcvf}. Figure~\ref{fig:rcv1sgdVsqn} shows that for each batch size,
the SQN method outperforms SGD, and both methods improve as batch size increases. We observe that using  $\bb=300$ or 1000 yields different relative outcomes for the SQN method  when measured in terms of iterations or adp: a batch size of 300 provides the fastest initial decrease in the objective, but that method is eventually overtaken by the variant with the larger batch size of 1000. 

\begin{figure}[H]
	\centering
	\subfloat[Varying $\bbh$ for $\bb=300$]
	{ \includegraphics[scale=.60]{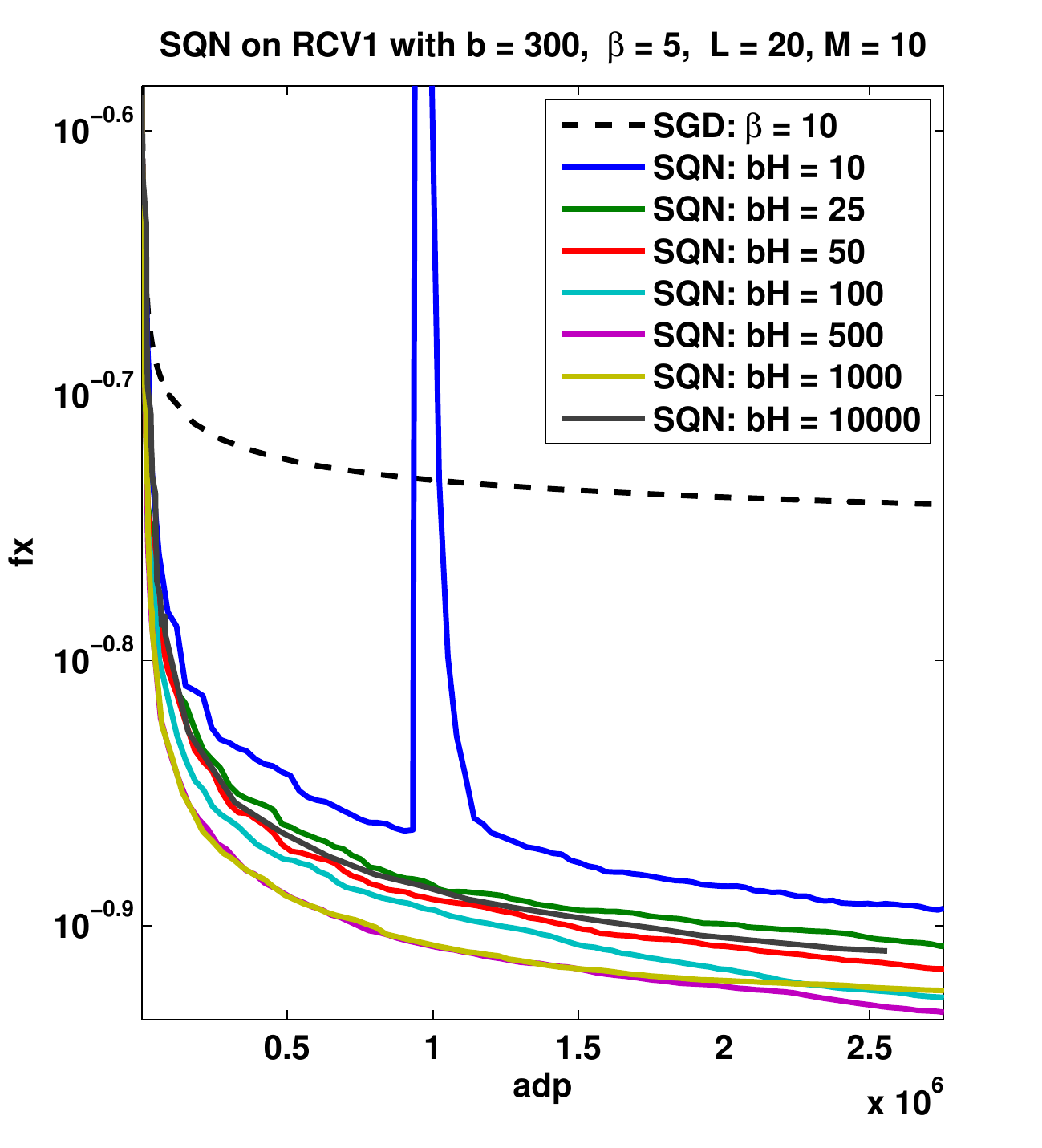}}
	\subfloat[Varying $\bbh$ for $\bb=1000$]
	{\hbox{\hspace{-0.95ex}\includegraphics[scale=.60]{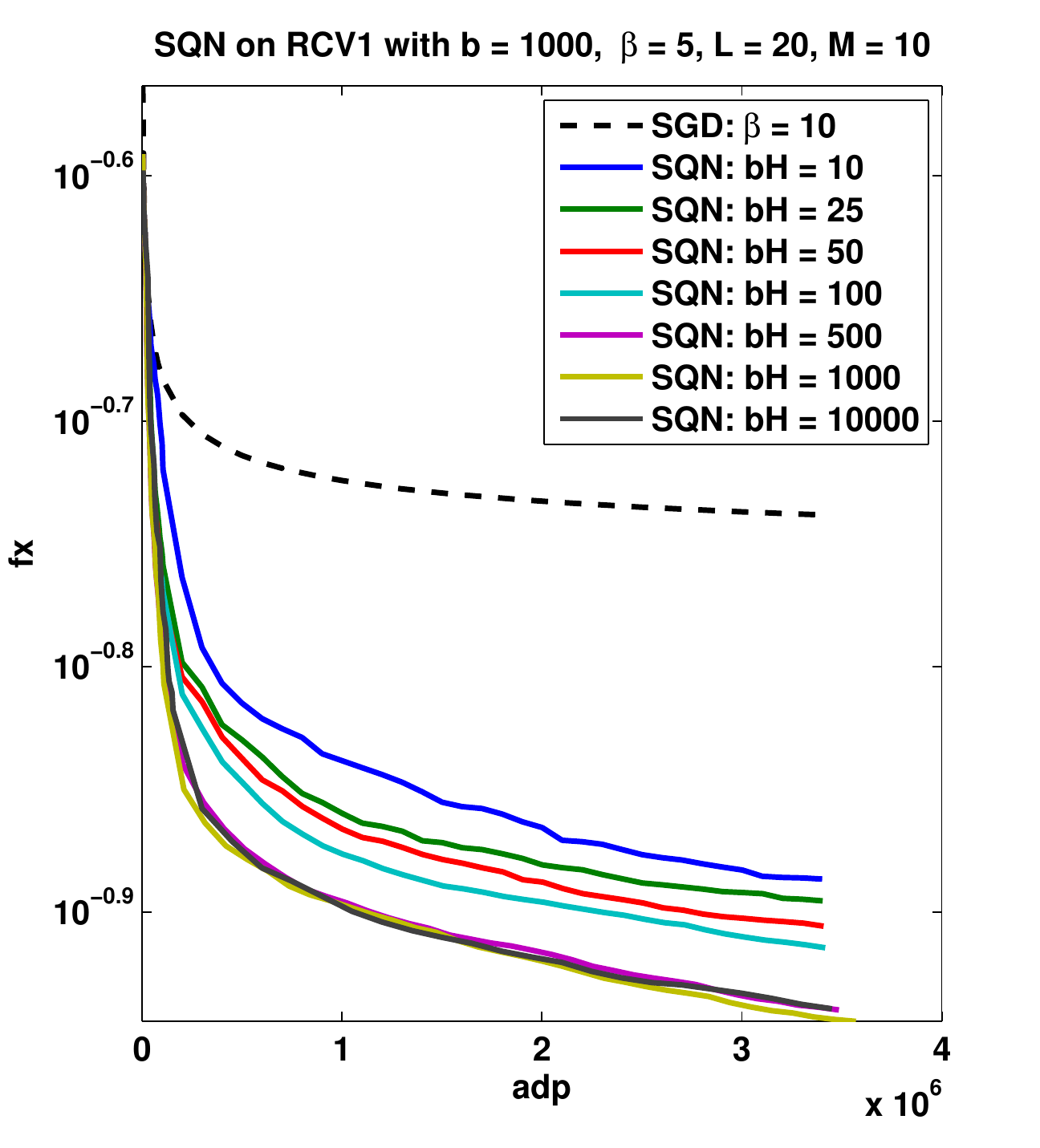}}}
		\caption{Varying Hessian batch size parameter $\bbh$ on the RCV1 dataset for gradient batch values $\bb$ of 300 and 1000. All other parameters in the SQN method are held constant
		at $\bL=20$, $\bM=10$, and $\beta=5$. }
	\label{fig:RCV1varyBHV}
\end{figure}

Figure~\ref{fig:RCV1varyBHV} illustrates the effect of varying the Hessian batch size $\bbh$ from 10 to 10000, while keeping the gradient batch size $\bb$ fixed at 300 or 1000. For $\bb=300$ (Figure~\ref{fig:RCV1varyBHV}a) increasing $\bbh$ improves the performance of SQN, in terms of adp, up until $\bbh=1000$, where the benefits
of the additional accuracy in the Hessian approximation do not outweigh the additional computational cost.
In contrast, Figure~\ref{fig:RCV1varyBHV}b shows that for $\bb=1000$, a high value for $\bbh$, such as 10000, can be effectively used since the cost of the Hessian-vector product relative to the gradient is lower. One of the conclusions drawn from this experiment is that there is much freedom in the choice of $\bbh$, and that only a small sub-sample  of the data (e.g. $\bbh = 100$) is needed for the stochastic quasi-Newton approach to yield benefits.

One should guard, however,  against the use of very small values for $\bbh$, as seen in the large blue spike in Figure~\ref{fig:RCV1varyBHV}a corresponding to $\bbh=10$.
To understand this behavior, we monitored the angle between the vectors $s$ and $y$ and observed that it approached $90^\circ$ 
between iteration 3100 and 3200, which is where the spike occurred.
Since the term $ s^T y$ enters in the denominator of the BFGS update formula \eqref{lbfgs2}, this led to a very large and poor step. Monitoring $ s^T y$ (relative to, say, $s^T B s$)  can a useful indicator of a potentially harmful update; one can increase $ \bbh$ or skip the update when this number is smaller than a given threshold.

\begin{figure}[H]
\centering
\includegraphics[scale=.60]{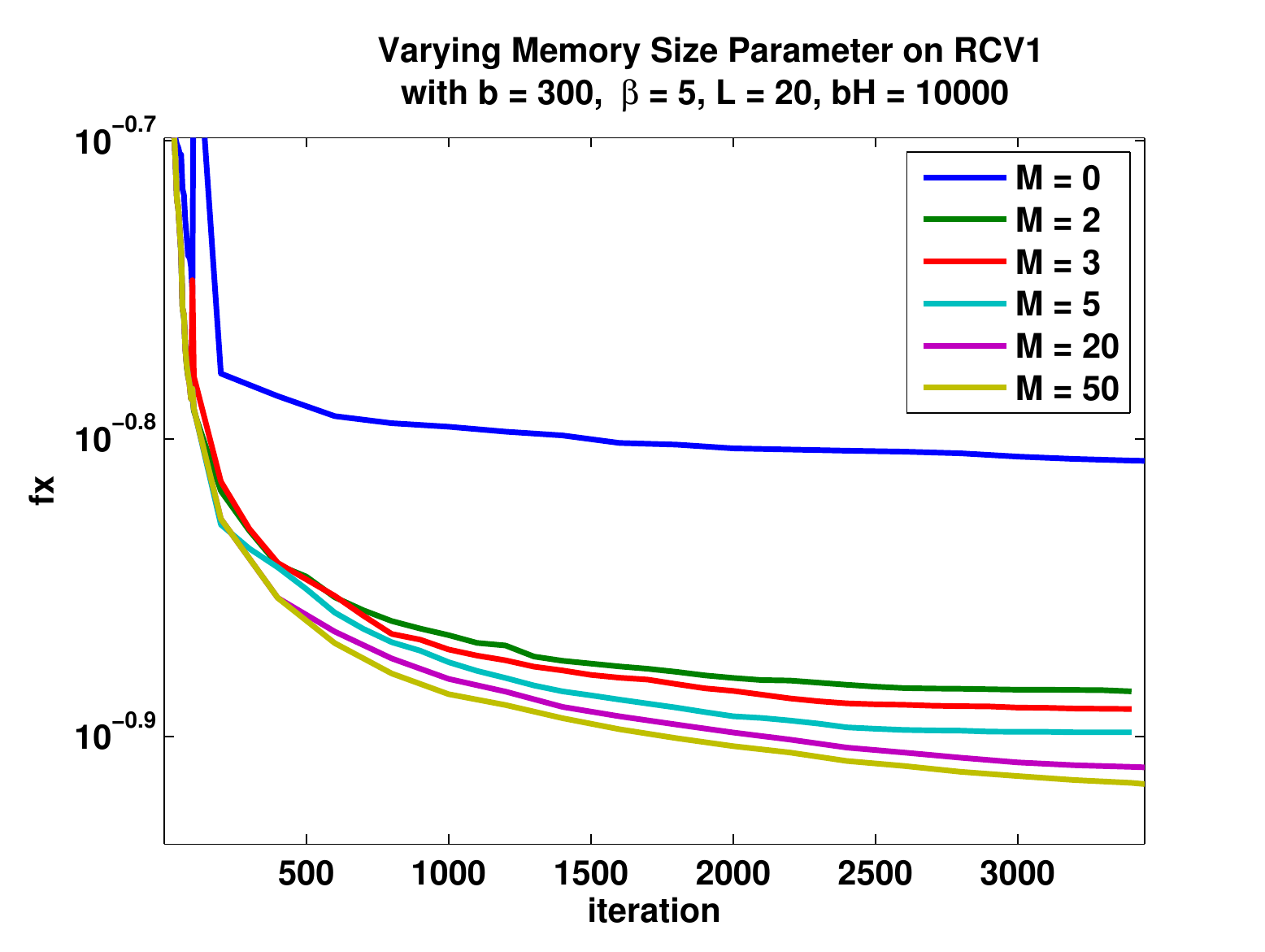}
\caption{Impact of the memory size parameter on the RCV1 dataset. $\bM$ is varied between 0 and 50 while all other parameters are held constant at
$\bb= 300$, $\bL=20$, and $\bbh=10000$.}
\label{fig:rcv1varyM_b300}
\end{figure}

The impact of the memory size parameter $\bM$ is shown in Figure~\ref{fig:rcv1varyM_b300}. The results improve consistently as $\bM$ increases, but beyond $\bM=2$ these improvements are rather small,
especially in comparison to the results in Figure \ref{random2} for the synthetic data. 
The reasons for this difference are not clear, but for the deterministic L-BFGS method the effect of $\bM$
on performance is known to be problem dependent. We observe that performance with a value of $\bM=0$, which results in a Hessian approximation of the form
$H_t = \frac{s_t^Ty_t}{y_t^Ty_t}I$,
is poor and also unstable in early iterations, as shown by the spikes in Figure~\ref{fig:rcv1varyM_b300}.

To gain a better understanding of the behavior of the SQN method, we also monitored the following two errors:
\begin{equation}
 \label{rel1}
\textrm{GradError} = \frac { \left\| \nabla F(w) - \nh  F(w) \right\|_2}{ \left\| \nabla F(w) \right\|_2},   
\end{equation}
and
\begin{equation}
\label{rel2}
\textrm{HvError} = \frac { \left\| \nabla^2 F(\bar{w}_I)(\bar{w}_I-\bar{w}_J) - \nh^2 F(\bar{w}_I)(\bar{w}_I-\bar{w}_J) \right\|_2}{ \left\| \nabla^2 F(\bar{w}_I)(\bar{w}_I-\bar{w}_J) \right\|_2} .
\end{equation}

\medskip\noindent The quantities $\nabla F(w)$ and $\nabla^2 F(\bar{w}_I)(\bar{w}_I-\bar{w}_J)$ are computed with the entire data set, as indicated by \eqref{rcvf}. Therefore, the ratios above report the relative error in the stochastic gradient used in \eqref{stqn} and the relative error in the computation of the Hessian-vector product \eqref{newcor2}.

\begin{figure}[H]
 \centering
\includegraphics[scale=.65]{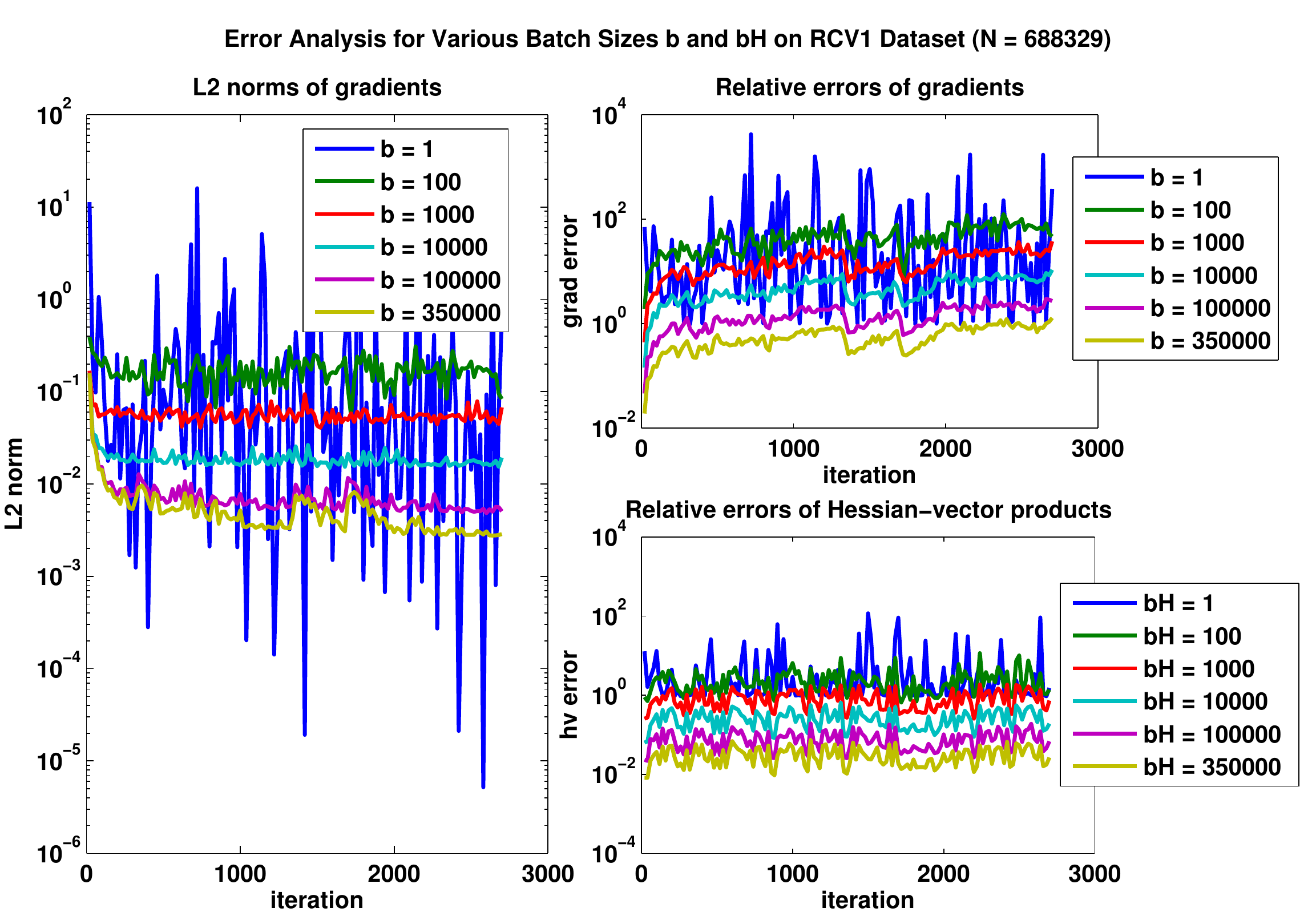}
\caption{Error plots for RCV1 dataset. The figure on the left plots $\| \nh F(w) \|_2$ for various values of $\bb$. The figures on the right display the errors \eqref{rel1} and \eqref{rel2}.
The errors were calculated every 20 iterations during a single run of SQN with parameters $\bb=300$, $\bL=20$, $\bM=5$, and $\bbh= 688329$.}
\label{rcv1:errors1}
\end{figure}

Figure~\ref{rcv1:errors1} displays these relative errors for various batch sizes $\bb$ and $\bbh$, along with the norms of the stochastic gradients. These errors were calculated every 20 iterations during a \emph{single run} of SQN with the following parameter settings:  $\bb=300$, $\bL=20$, $\bM=5$, and $\bbh=688329$.
Batch sizes larger than $\bb=10000$ exhibit non-stochastic behavior in the sense that all relative errors are less than one, and the norm of these approximate gradients decreases during the course of the iteration. Gradients with a batch size less than 10000 have relative errors greater than 1, and their norm does not exhibit
decrease over the course of the run.

The leftmost figure also shows that the $\ell_2$ norms of the stochastic gradients
decrease as the batch size $\bb$ increases, i.e., there is a tendency for inaccurate gradients to have a larger norm,  as expected from the geometry of the error.


Figure~\ref{rcv1:errors1} indicates that the Hessian-vector errors stay relatively constant throughout the run and have smaller relative error than the gradient.  As discussed above, some accuracy here is important while it is not needed for the batch gradient.

\subsection{A Speech Recognition Problem} \label{speechs}

The speech dataset, provided by Google, is a collection of feature vectors representing 10 millisecond frames of speech with a corresponding label representing the phonetic state assigned to that frame. Each feature $x_i$ has a dimension of {\it NF} $= 235$  
and has corresponding label $z_i\in C= \left\{1, 2, \ldots, 129 \right\}$. There are a total of $N=191607$ samples; the number of variables is $n=$ {\it NF} $\times |C|= 30315$.

The problem is modeled using multi-class logistic regression. The unknown parameters are assembled in a matrix $W\in\mathbb{R}^{|C| \times \textit{NF}}$, and the objective
is given by
\begin{align}
F(W) &= -\frac{1}{N}\sum_{i=1}^N \log\left( \frac{\exp( W_{z_i}x_i)}{\sum_{j\in\mathcal{C}}\exp(W_jx_i)}\right),  \label{multi}
\end{align}
where $x_i\in\mathbb{R}^{\textit{NF} \times 1}$, $z_i$ is the index of the correct class label, and 
$W_{z_i} \in\mathbb{R}^{1\times \textit{NF}}$ is the row vector corresponding to the weights associated with class $z_i$.

\begin{figure}[H]
\centering
\includegraphics[scale=.65]{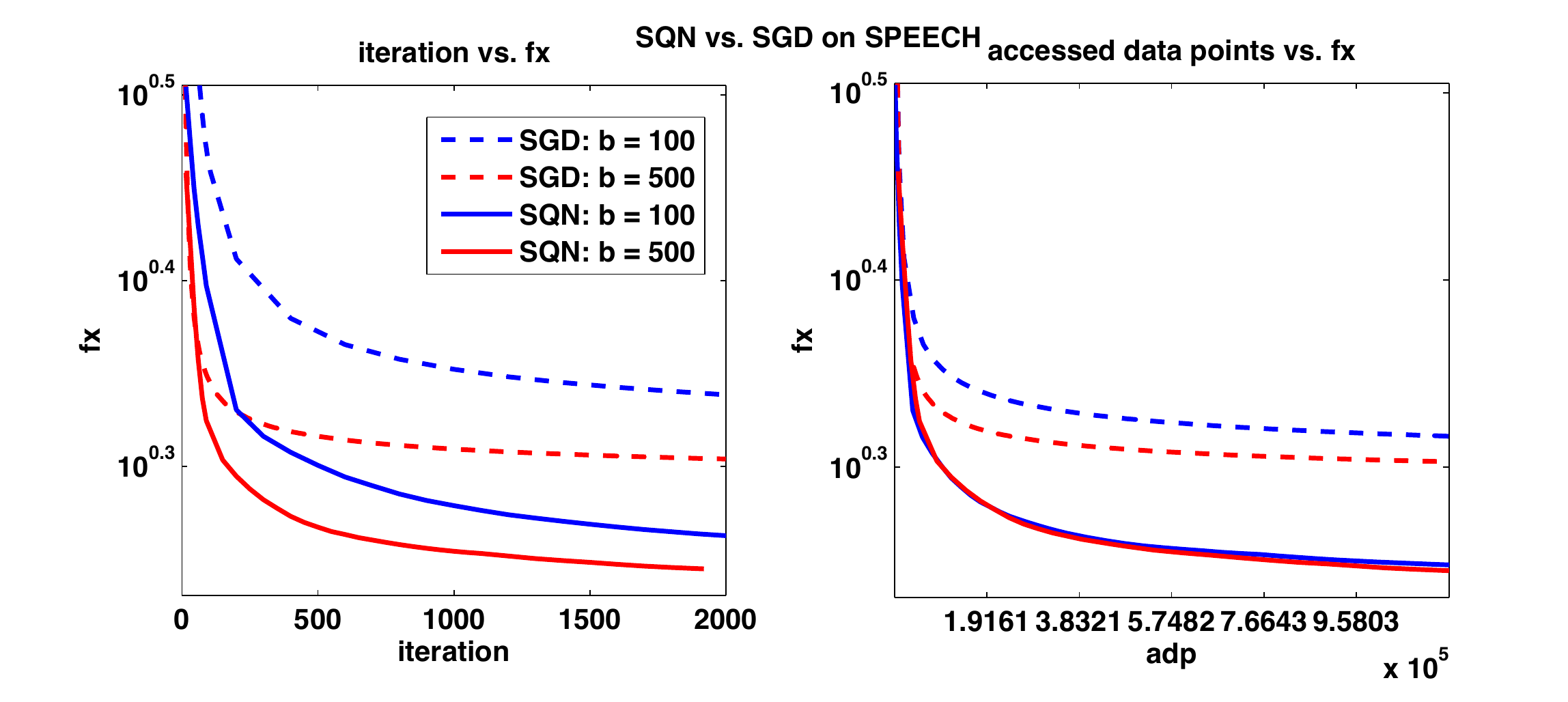}
\caption{Illustration of SQN and SGD on the SPEECH dataset. The gradient batch size $\bb$ is  100 or 500, and for SQN we use $\bbh=1000$, $\bM=5$ and $\bL=10$. 
In the rightmost graph, the tick marks on the x-axis denote the epochs of the SGD method. }
\label{fig:speechsgdVsqn}
\end{figure}

\medskip
Figure \ref{fig:speechsgdVsqn} displays the performance of SGD and SQN for $\bb= 100$ and 500 (which represent approximately 0.05\%, and  0.25\% of the dataset).  For the SGD method, we chose the step length $\beta=10$ for both values of $\bb$; for the SQN method  we set $\beta=2$,
$\bL=10, \ \bM=5, \  \bbh=1000$.

We observe from Figure~\ref{fig:speechsgdVsqn} that SQN improves upon SGD in terms of adp, both initially and in finding a lower objective value. 
Although the number of SQN iterations decreases when $\bb$ is increased from 100 to 500,  in terms of computational cost the two versions of the SQN method yield almost identical performance.
%

\begin{figure}
\centering
\includegraphics[scale=.60]{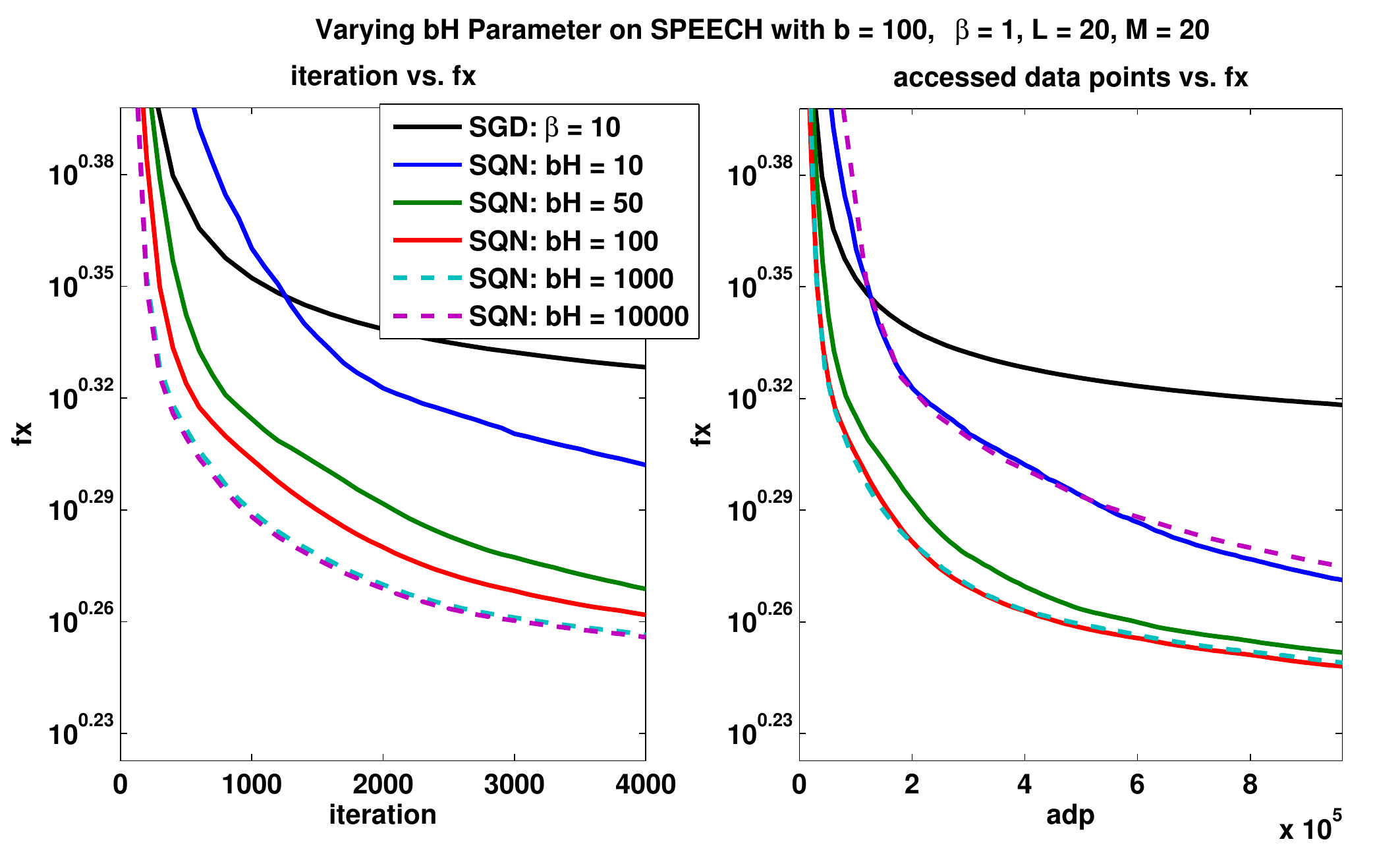}
\caption{Varying Hessian batch size parameter $\bbh$ on SPEECH dataset. All other parameters
are held constant at $\bb=100$, $\bL=20$, $\bM=20$.}
\label{speechVaryBHV}
\end{figure}

The effect of varying the Hessian batch size $ \bbh$  is illustrated in Figure~\ref{speechVaryBHV}. The figure on the left shows that increasing $ \bbh$ improves the performance 
of SQN, as measured by iterations, but only marginally from $ \bbh=1000$ to $ 10,000$. Once the additional computation cost of Hessian-vector products is accounted for, we observe from the figure on the right that $\bbh =100$ is as effective as $\bbh=1000$.
Once more, we conclude that only a small subset of data points $\Sh$ is needed to obtain useful curvature information in the SQN method.


\begin{SCfigure}
 \centering
\includegraphics[scale=.55]{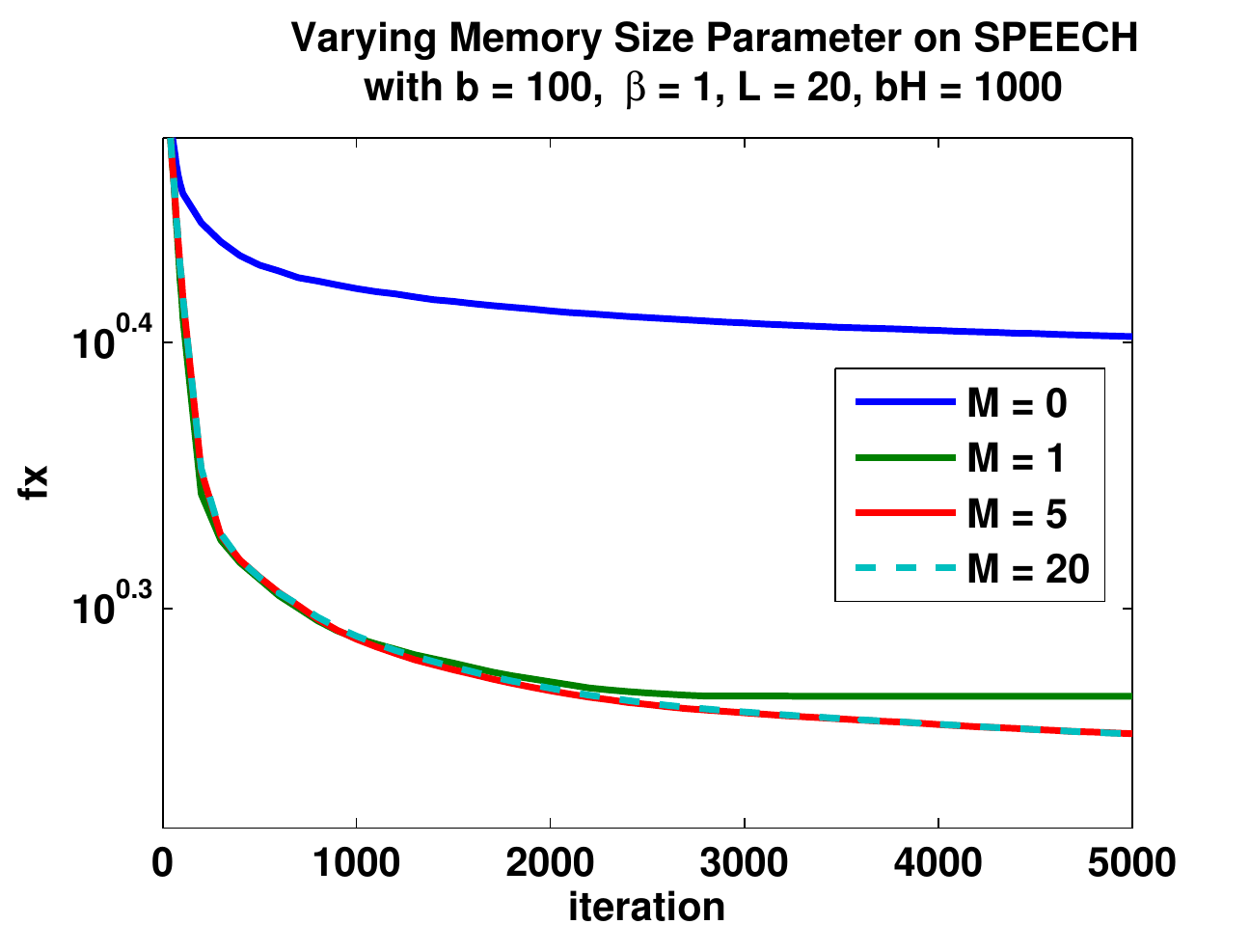}
\caption{Performance on the SPEECH dataset with varying memory size $\bM$. All other parameters
are held constant at $\bb=100$, $\bL=20$, $\bbh=1000$.  \\ \\ }
\label{speechVaryM}
\end{SCfigure}

Figure \ref{speechVaryM} illustrates the impact of increasing the memory size $\bM$ from 0 to 20 for the SQN method.
A memory size of zero leads to a marked degradation of performance. Increasing $\bM$ from 0 to 5 improves SQN, but values greater than 5 yield no measurable benefit.


\subsection{Generalization Error}

\begin{figure}[H]
	\centering
	\subfloat[Training and test error on the objective]
	{ \includegraphics[scale=.55]{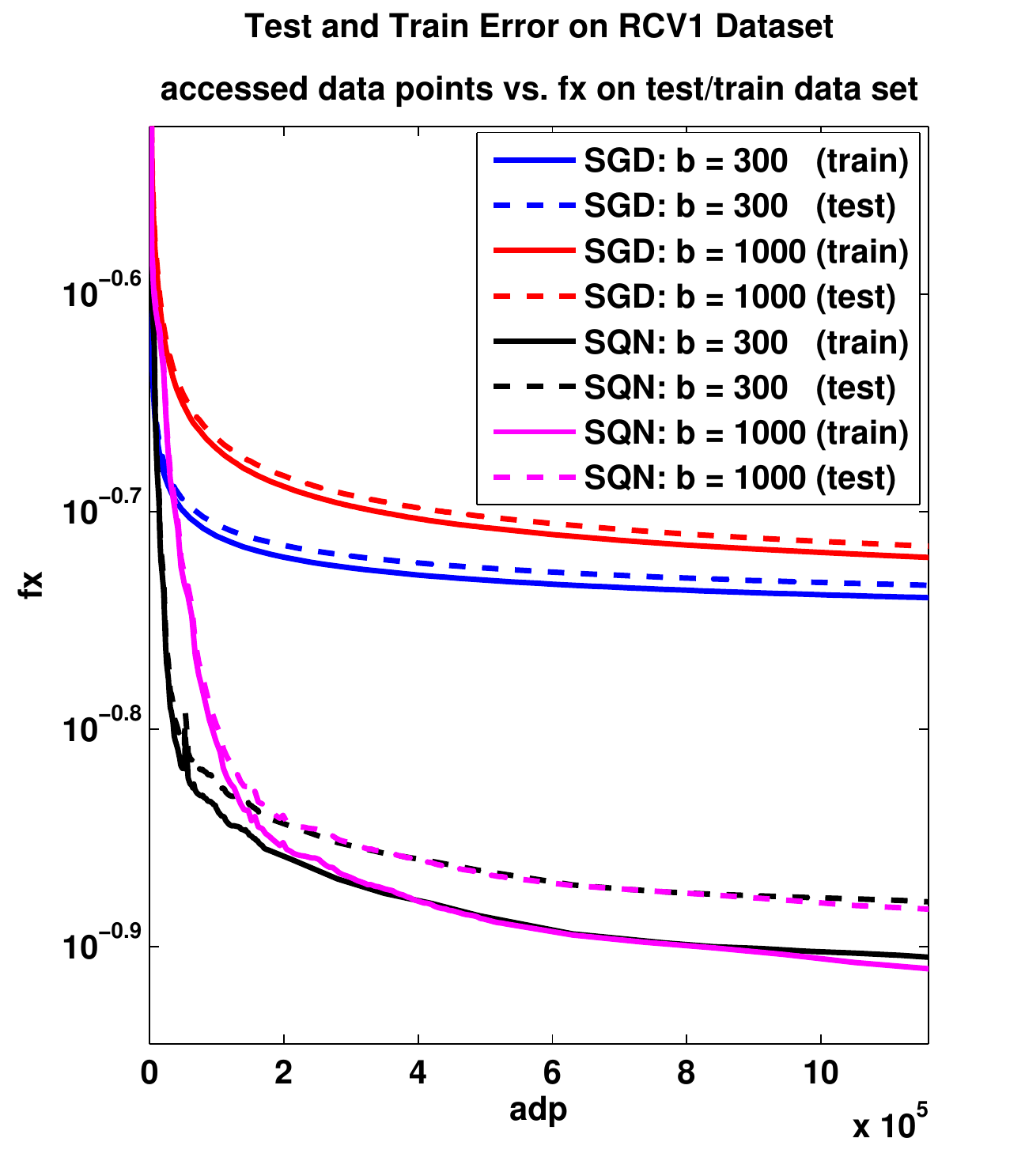}}
	\subfloat[Test error for percent correctly classified]
	{\hbox{\hspace{-0.95ex}\includegraphics[scale=.55]{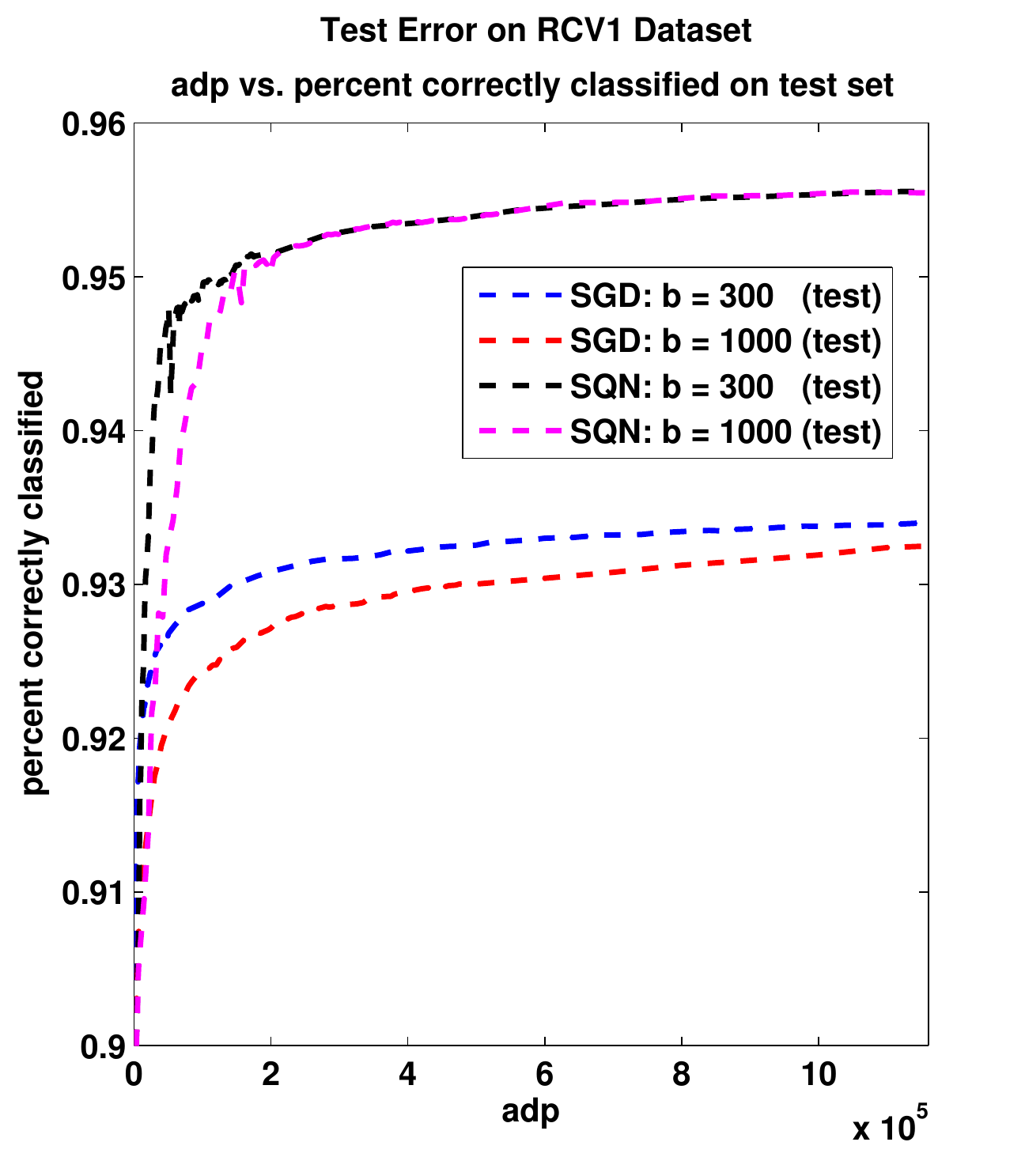}}}
		\caption{Illustration of the generalization error on the RCV1 dataset. For both SGD and SQN, $\bb$ is set to 300 or 1000; for SQN we set $\bbh=1000$,  $\bM=5$ and $\bL=20$.  }
	\label{fig:rcv1_trainingError}
\end{figure}

The primary focus of this paper is on the minimization of training error \eqref{empirical}, but it is also interesting to explore the performance of the SQN method in terms of generalization (testing) error. For this purpose we consider the RCV1 dataset, and in
Figure~\ref{fig:rcv1_trainingError} we report the performance of algorithms SQN and SGD with respect to unseen data (dotted lines). Both algorithms
were trained using $75\%$ of the data and then tested on the remaining $25\%$ (the test set). In Figure~\ref{fig:rcv1_trainingError}a,
the generalization error is measured in terms of decrease of the objective (\ref{rcvf}) over the test set, and in Figure~\ref{fig:rcv1_trainingError}b, in terms of the percent of correctly classified data points from the test set. The first measure complements the latter in the sense that it takes into account the confidence of the correct predictions and the inaccuracies wrought by the misclassifications. Recall that there are 2 classes in our RCV1 experiments, so random guessing would yield a percentage of correct classification of 0.5.

As expected, the objective on the training set is lower than the objective on the test set, but not by much. These graphs suggests that over-fitting is not occurring since the objective on the test set decreases monotonically. The performance of the stochastic quasi-Newton method is clearly very good on this problem.

\subsection{Small Mini-Batches}

In the experiments reported in the sections ~\ref{rcv1s} and \ref{speechs}, we used fairly large gradient batch sizes,  such as $\bb=50, 100, 1000$, because they gave good performance for both the SGD and SQN methods on our test problems. Since we set $\bM=5$, the cost of the multiplication $H_t \widehat \nabla F(w^k)$ (namely, $4 \bM n= 20n$)  is small compared to the cost of $\bb n$ for a batch gradient evaluation. We now explore the efficiency of the stochastic quasi-Newton method for smaller values of the batch size $\bb$. 

\begin{figure}[H] 
 \centering
\includegraphics[scale=.70]{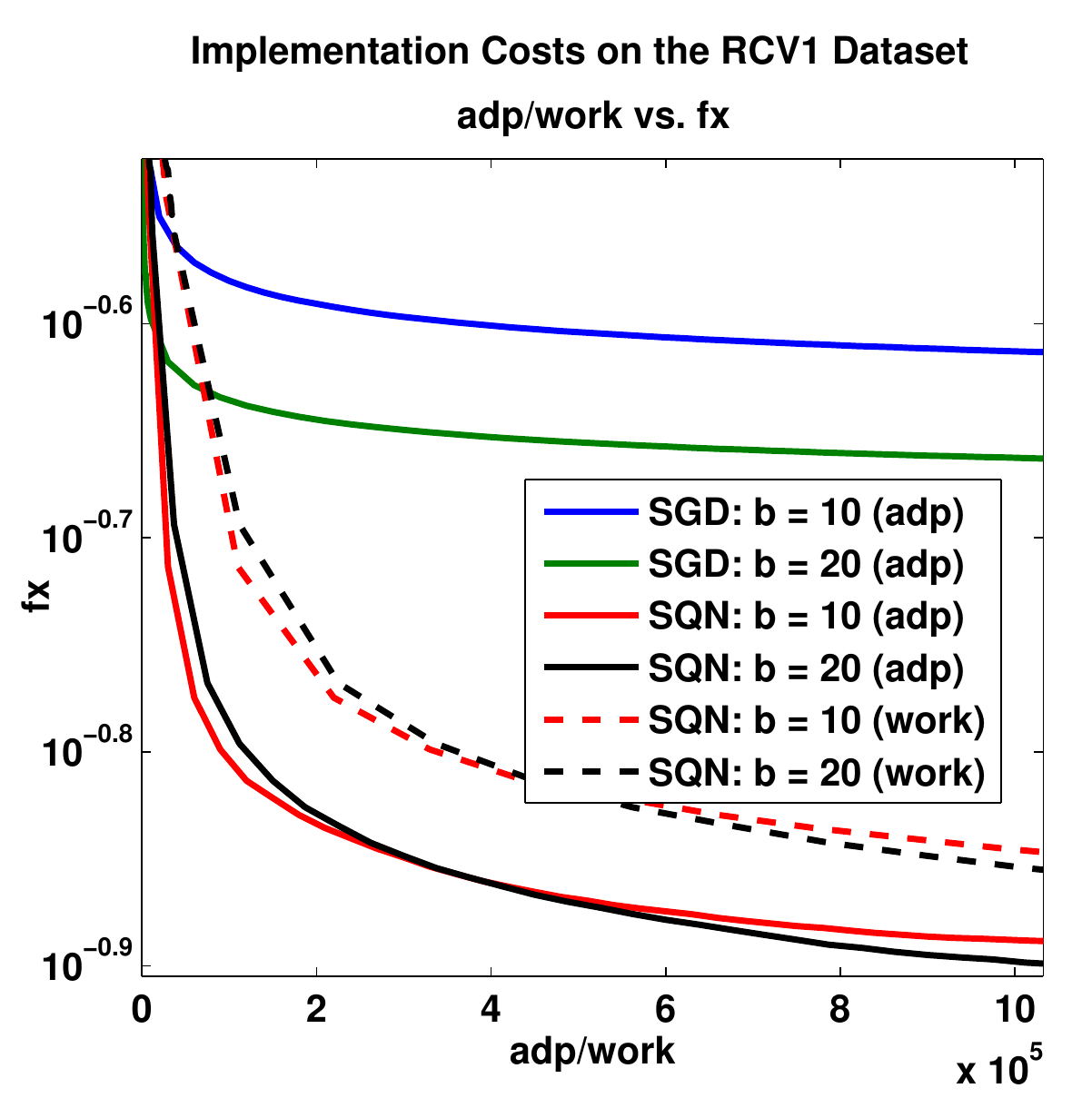}
\caption{Comparison using $\bb=10, \ 20$ on RCV1. The solid lines measure performance in terms of adp and the dotted lines measure performance in terms of total computational work \eqref{counts} (scaled by a factor of $1/n$).
For SQN we set $\bM=5$, $\bbh=1000$, $\bL=200$, $\beta = 1$, and
for SGD we set $\beta = 5$. 
}
\label{bsmall}
\end{figure}

In Figure~\ref{bsmall} we report results for the SGD and SQN methods for problem RCV1, for $\bb = 10$ and 20. We use two measures of performance: total \emph{computational work} and adp.
For the SQN method, the work measure is given by \eqref{counts}, which includes
the evaluation of the gradient \eqref{bat}, the computation of the quasi-Newton step \eqref{stqn}, and the Hessian-vector products \eqref{newcor2}.

 
In order to compare total work and adp on the same figure, we scale the work by $1/n$.
The solid lines in Figure~\ref{bsmall} plot the objective value vs adp, while the dotted lines plot function value vs total work. 
We observe from Figure~\ref{bsmall} that,  even for small $\bb$,
the SQN method outperforms SGD by a significant margin  in spite of the additional  Hessian-vector product cost.
Note that in this experiment the $4Mn$ cost of computing the steps is still less than half the total computational cost \eqref{counts} of the SQN iteration.

In this experiment, it was crucial to update  the quasi-Newton matrix infrequently ($\bL= 200$), as this allowed us to employ a large value of $ \bbh$ at an acceptable cost. In general, the parameters $\bL, \ \bM$ and $ \bbh$ provide much freedom in adapting the SQN method to a specific application.

%
%
%
%
\subsection{Comparison to the oLBFGS method}
We also compared our algorithm to the oLBFGS method  \cite{schraudolph2007stochastic}, which is the best known stochastic quasi-Newton method in the literature. It is of the form \eqref{direct} but differs from our approach in three crucial respects: the L-BFGS update is performed at every iteration,  the curvature estimate is calculated using gradient differencing, and the sample size for gradient differencing is the same as the sample size for the stochastic gradient. This approach requires  two gradient evaluations per iteration; it computes
\[ 
w^{k+1} = w^k - \alpha^k H_k \nh F_{\mathcal{S}_k}(w^k), \ \ s_k = w^{k} - w^{k-1}, \  \ y_k = \nh F_{\mathcal{S}_{k-1}}(w^k)-  \nh F_{\mathcal{S}_{k-1}}(w^{k-1}),
\]
where we have used subscripts to indicate the sample used in the computation of the gradient $\hat \nabla F$. The extra gradient evaluation is  similar in cost to our Hessian-vector product, but we compute that product only every $L$ iterations. Thus, the oLBFGS method is analogous to our algorithm with $L=1$ and $\bb= \bbh$, which as the numerical results below show, is not an efficient allocation of effort. In addition, the oLBFGS method is limited in the choice of samples $\cal S$ because, when these are small, the Hessian approximations may be of poor quality.

%

We implemented the oLBFGS method as described in \cite{schraudolph2007stochastic}, with the following parameter settings: i) we found it to be unnecessary to add a damping parameter  to the computation $y_k$, and thus set $\lambda=0$ in the reset $y_k \leftarrow y_k + \lambda s_k$; ii) the  parameter $\epsilon$ used to rescaled the first iteration, $w^1 = w^0 - \epsilon\alpha^k\nh F(w^0)$, was set to $\epsilon = 10^{-6}$; iii) the initial choice of scaling parameter in Hessian updating (see Step 1 of Algorithm~\ref{alg2}) was the average of the quotients $s_i^Ty_i/y_i^T y_i$ averaged over the last $M$ iterations, as recommended in \cite{schraudolph2007stochastic}.

\begin{figure}[H] 
 \centering
\includegraphics[scale=.65]{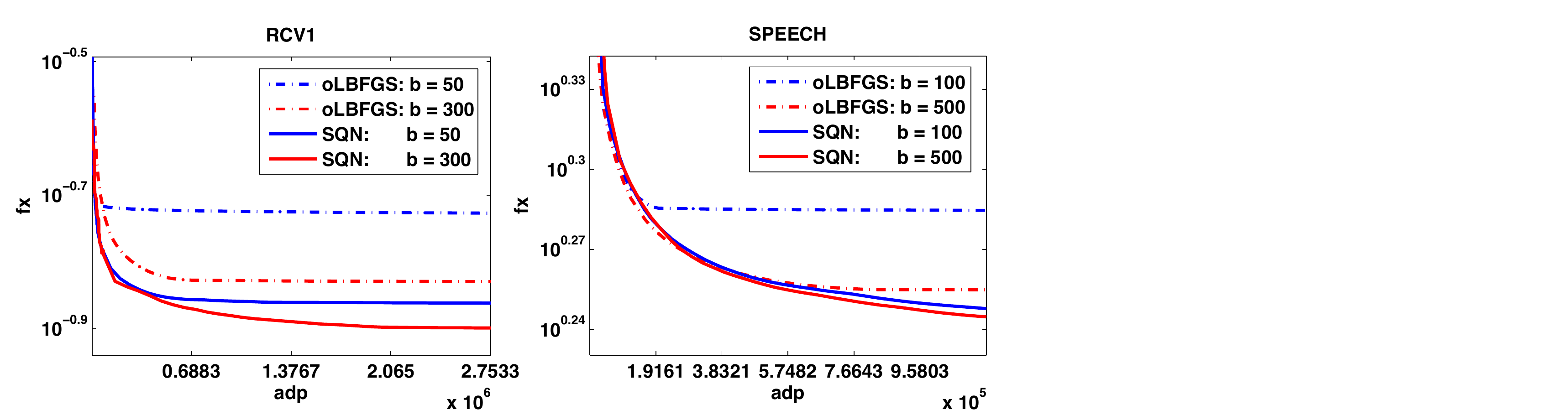}
  \caption{Comparison of oLBFGS (dashed lines) and SQN (solid lines) in terms of accessed data points. For RCV1 dataset gradient batches are set to $\bb = 50$ or $300$, for both methods; additional parameter settings for SQN are $L=20$, $\bbh = 1000$, $M=10$. For Speech dataset we set to $\bb = 100$ or $500$; and for SQN we set  $L=10$,  $\bbh = 1000$, $M=10$.}
  \label{fig:olbfgs}
\end{figure}

Figure~\ref{fig:olbfgs} compares our SQN method to the aforementioned oLBFGS on our two realistic test problems, in terms of accessed data points. We observe that SQN has overall better performance, which is more pronounced for smaller batch sizes.  

\section{Related Work}    \label{related}

Various stochastic quasi-Newton algorithms have been proposed in the literature \cite{schraudolph2007stochastic,mokhtariregularized,bordes2009sgd,roux2010fast}, 
but have not been entirely successful. The  methods in \cite{schraudolph2007stochastic} and \cite{mokhtariregularized} use the BFGS framework; the
first employs an L-BFGS implementation,  as mentioned in the previous section, and the latter uses a regularized BFGS matrix.
Both methods enforce uniformity in gradient differencing by resampling data points
so that two consecutive gradients are evaluated with the same sample $\Ss$; this strategy requires an extra gradient evaluation at each iteration. 
The algorithm
presented in \cite{bordes2009sgd} uses SGD with a diagonal rescaling matrix based 
on the secant condition associated with quasi-Newton methods. Similar to our approach, \cite{bordes2009sgd} updates the rescaling matrix
at fixed intervals in order to reduce computational costs.
A common feature of \cite{schraudolph2007stochastic,mokhtariregularized,bordes2009sgd} is that 
the Hessian approximation might be updated with a  high level of noise.

A two-stage online Newton strategy is proposed in \cite{bach2013nonT}. The first stage runs averaged SGD with a step size of order $O(1/\sqrt{k})$, and the second stage minimizes a quadratic model of the objective function using SGD with a constant step size. The second stage effectively takes one Newton step, and employs Hessian-vector products in order to compute stochastic derivatives of the quadratic model. This method is significantly different from our quasi-Newton approach.

A stochastic approximation method that has shown to be effective in practice is AdaGrad \cite{duchihazansing}. The iteration  is of the form \eqref{rm}, where $B_k$ is a diagonal matrix that estimates the diagonal of the squared root of the uncentered covariance matrix of the gradients; it is shown in \cite{duchihazansing} that such a matrix  minimizes a regret bound. The algorithm presented in this paper is different in nature from AdaGrad, in that it employs a full (non-diagonal) approximation to the Hessian $\nabla^2 F(w)$. 

Amari \cite{amari1998natural} popularized the idea of incorporating information from the geometric space of the inputs into online learning with his presentation of the natural gradient method. This method seeks to find the steepest descent direction in the feature space $x$ by using the Fisher information matrix, and is shown to achieve asymptotically optimal bounds. The method does, however, require knowledge of the underlying distribution of the training points $(x, z)$, and the Fisher information matrix must be inverted. These concerns are addressed in \cite{park2000adaptive}, which presents  an adaptive method for computing the inverse Fisher Information matrix in the context of multi-layer neural networks.

The authors of TONGA \cite{roux2007topmoumoute} interpret natural gradient descent as the direction that maximizes the probability of reducing the generalization error. They outline an online implementation using the uncentered covariance matrix of the empirical 
gradients that is updated in a weighted manner at each iteration. Additionally, they show how to maintain a low rank approximation of the covariance matrix so
that the cost of the natural gradient step is $O(n)$. 
In \cite{roux2010fast} it is argued that an algorithm should contain information about both the Hessian and covariance matrix, maintaining that
that covariance information is needed to cope with the variance due to the space of inputs, and Hessian information is useful to improve the optimization.

Our algorithm may appear at first sight to be similar to the method proposed by Byrd et al. \cite{ByrdChinNeveiNoce,dss}, which also employs Hessian-vector products to gain curvature information. We note, however, that the algorithms are different in nature, as the algorithm presented here operates in the stochastic approximation regime, whereas \cite{ByrdChinNeveiNoce,dss} is a batch (or SAA) method. 

\section{Final Remarks}   \label{final}
In this paper, we presented a quasi-Newton method that operates in the stochastic approximation regime. It is designed for the minimization of convex stochastic functions, and was tested on problems arising in machine learning. In contrast to previous attempts at designing stochastic quasi-Newton methods, our approach does not  compute gradient differences at every iteration to gather curvature information; instead it computes (sub-sampled) Hessian-vector products at regular intervals to obtain this information in a stable  manner.

Our numerical results suggest that the method does more than rescale the gradient, i.e., that  its improved performance over the stochastic gradient descent method of Robbins-Monro is the result of incorporating curvature information in the form of a full matrix.

The practical success of the algorithm relies on the fact that the batch size $\bbh $ for Hessian-vector products can be chosen large enough to provide useful curvature estimates, while the update spacing $L$ can be chosen large enough (say $L= 20$) to amortize the cost of Hessian-vector products, and make them affordable.
Similarly, there is a wide range of values for the gradient batch size $b$ that makes the overall quasi-Newton approach (\ref{direct}) viable.
 
The use of the Hessian-vector products \eqref{newcor2} may not be essential; one might be able to achieve the same goals using differences in gradients, i.e.,
\[
\bar y = \bg(\bar w_t ) - \bg (\bar w_{t-1}).
\]
This would require, however, that the evaluation of these gradients employ the same sample $\cal S$ so as to obtain sample uniformity, as well as the development of a strategy to prevent close gradient differences from magnifying round off noise. In comparison, the use of Hessian-vector products takes care of these issues automatically, but it requires code for Hessian-vector computations (a task that is not often not onerous). 

We established global convergence of the algorithm on strongly convex objective functions.  Our numerical results  indicate that the algorithm is more effective than the best known stochastic quasi-Newton method (oLBFGS \cite{schraudolph2007stochastic}) and
suggest that it holds much promise for the solution of large-scale problems arising in stochastic optimization. Although we presented and analyzed the algorithm in the convex case, our approach   is applicable to non-convex problems provided it employs a mechanism for ensuring that the  condition $s_t^T y_t >0$ is satisfied.


\newpage
\small
\bibliographystyle{plain}

\bibliography{../../References/references}
\end{document}